\documentclass[twoside,11pt]{article}
\usepackage{jmlr2e}

\usepackage{amsmath}
\usepackage{mathtools}
\usepackage{enumitem}
\usepackage{placeins}
\usepackage{xcolor}
\usepackage{tabularx}
\usepackage{multirow}
\usepackage{booktabs}
\usepackage[compatibility=false]{caption}
\usepackage{pifont}
\usepackage{algorithm}
\usepackage{algorithmic}
\usepackage{tikz}
\usepackage[utf8]{inputenc}
\usepackage{hyperref}

\usepackage{ragged2e}
\usepackage[american]{babel}
\usepackage{subcaption}
\DeclareCaptionLabelFormat{custom}{(#2)}
\newcommand{\E}{\mathbb{E}}
\AtBeginDocument{\renewcommand{\Pr}{\mathbb{P}}}
\newcommand{\R}{\mathbb{R}}
\usepackage{amssymb}

\newcommand{\cX}{\mathcal{X}}
\newcommand{\cY}{\mathcal{Y}}
\newcommand{\cF}{\mathcal{F}}

\newcommand{\fstar}{f^\star}

\newcommand{\KL}{\mathrm{KL}}
\newcommand{\TV}{\mathrm{TV}}
\newcommand{\cN}{\mathcal{N}}

\let\proof\relax

\makeatletter
\let\c@theorem\relax    %
\makeatother
\usepackage{amsthm}
\newtheorem{theorem}{Theorem}[section]
\newtheorem{lemma}[theorem]{Lemma}
\newtheorem{definition}[theorem]{Definition}
\newtheorem{remark}[theorem]{Remark}

\ShortHeadings{Conformal Prediction at Performative Equilibrium}{Mahrous and Di Pietro}
\firstpageno{1}

\begin{document}

\title{Conformal Prediction at Performative Equilibrium: The Price of Performativity}

\author{\name Ahmed Mahrous \email ahmed.mahrous@kaust.edu.sa \\
       \addr King Abdullah University of Science and Technology (KAUST)\\
       Thuwal 23955, Saudi Arabia
       \AND
       \name Roberto Di Pietro \email roberto.dipietro@kaust.edu.sa \\
       \addr King Abdullah University of Science and Technology (KAUST)\\
       Thuwal 23955, Saudi Arabia}

\editor{} %

\maketitle

\begin{abstract}
Conformal prediction gives finite-sample coverage when calibration and test points are exchangeable, but exchangeability can fail when deployment changes the data distribution. We study conformal prediction in performative settings, where the deployed predictor affects the data distribution, and seek prediction sets valid at performative equilibrium: coverage must hold under $D(\fstar)$ for a performatively stable predictor $\fstar$. Under total-variation $L$-sensitivity of the distribution map, we prove a phase transition at $L=1$. If $L \ge 1$, no finite-sample procedure can uniformly guarantee both equilibrium-valid coverage and nontrivial prediction sets. If $L<1$, the minimax excess miscoverage over informative procedures is $\widetilde{\Omega}(1/((1-L)\sqrt{n}))$, implying a sample-complexity inflation of order $1/(1-L)^2$ for fixed excess miscoverage. Our main technical tool is a transfer theorem: any procedure returning a predictor within $\varepsilon$ of equilibrium yields equilibrium-valid coverage up to an additive $L\varepsilon$ term after split-conformal calibration. We instantiate the theorem with two retraining algorithms: a model-free method for discretizable classes and a parametric tail-averaged method that matches the $n$- and $L$-dependence of this rate up to logarithmic factors. Experiments, including semi-synthetic credit scoring, demonstrate the predicted threshold behavior and scaling laws.\end{abstract}


\section{Introduction}
\label{sec:intro}

Conformal prediction provides distribution-free coverage guarantees:
given exchangeable calibration data, one can construct prediction sets
$C(x) \subseteq \cY$ satisfying
$\Pr(Y \in C(X)) \geq 1-\alpha$~\citep{vovk2005algorithmic,shafer2008tutorial,lei2018distributionfree,angelopoulos2023gentle}.  This guarantee requires no assumptions on the data-generating distribution
beyond exchangeability; informally, the test point must come from the same
data-generating process as the calibration data~\citep{vovk2012conditional,lei2014distribution}.

In many applications, this assumption fails because deployment changes
behavior and hence the data distribution.  A lending rule affects who receives credit and thus the future distribution
of defaults on which the rule is later retrained; a recidivism predictor
affects sentencing and thus the future distribution of reoffending; a
treatment-allocation rule affects who receives care and thus the future
distribution of health outcomes. In
such settings, the distribution is \emph{endogenous}: it depends on the
deployed predictor.  The performative prediction
framework~\citep{perdomo2020performative,mendler2020stochastic,miller2021echo,izzo2021performative,brown2022stateful,jagadeesan2022regret,hardt2022performative} formalizes this through a map
$D : \cF \to \Delta(\cX \times \cY)$ from predictors to distributions.
Deploying $f$ induces distribution $D(f)$, so a conformal guarantee
calibrated under a fixed distribution need not remain valid under
$D(f)$.

If $\fstar$ is a performatively stable predictor (a fixed point that minimizes risk under the distribution it induces), we want prediction sets valid at that distribution:
\begin{equation}\label{eq:equil-coverage}
\Pr_{(X,Y)\sim D(\fstar)}\!\bigl(Y \in C_{\fstar}(X)\bigr)
\;\geq\; 1 - \alpha.
\end{equation}
We call this \emph{equilibrium-valid coverage}.  Equilibrium is the
relevant operating point because deployed systems are retrained until
they stabilize, and coverage is needed at that stable distribution. This viewpoint aligns with prior work on retraining dynamics, optimization to stable points, performative risk minimization, and state-dependent feedback~\citep{perdomo2020performative,mendler2020stochastic,miller2021echo,izzo2021performative,brown2022stateful,jagadeesan2022regret}.

\paragraph{Research question.}
In performative environments, what coverage guarantees can conformal
prediction provide at performative equilibrium?  Under what conditions on the feedback map are such guarantees achievable, and what finite-sample complexity do they require?

\paragraph{Contributions.}
We answer these questions in five steps: we identify a feasibility
threshold, establish minimax lower and upper bounds, prove a transfer
theorem that converts convergence to equilibrium into coverage at
equilibrium, apply the theory with simple retraining-based
algorithms, and validate the predicted theory empirically.

\begin{enumerate}[leftmargin=*,itemsep=2pt]

\item \textbf{Impossibility (phase transition at $L\ge 1$). (Theorem~\ref{thm:impossibility}).}
If $L \ge 1$, no finite-sample procedure can guarantee \emph{nontrivial}
equilibrium-valid coverage over all $L$-sensitive performative environments:
there exists an environment for which the procedure must either fail to achieve
$1-\alpha$ coverage at equilibrium, or output the trivial prediction set that
includes all possible labels.

\item \textbf{Minimax rate (Theorem~\ref{thm:minimax}).}
When $L < 1$, the minimax excess miscoverage over informative procedures
(those whose prediction sets are $\tau$-informative, excluding the
vacuous set $C(x)=\cY$) with $n$ samples is
$\widetilde{\Omega}(1/((1{-}L)\sqrt{n}))$.  Under standard regularity ($\beta$-smooth, $\gamma$-strongly convex
loss with $\rho = L\beta/\gamma < 1$), Algorithm~\ref{alg:pcp}
achieves rate $\widetilde{O}(1/((1-\rho)\sqrt{n}))$, which matches
this rate when $\beta/\gamma = O(1)$. A model-free
algorithm achieves
$\widetilde{O}(1/((1{-}L)^{3/2}\sqrt{n}))$ under only $L < 1$.

\item \textbf{Coverage transfer theorem (Theorem~\ref{thm:meta}).}
If a learning or retraining procedure returns a predictor within
$\varepsilon$ of performative equilibrium in the metric used by the
$L$-sensitivity condition, then split conformal calibration on held-out
data achieves equilibrium-valid coverage up to an additive penalty of
order $L\varepsilon$ plus the standard finite-sample conformal term.
This decouples optimization from calibration, allowing equilibrium to be reached first and coverage to be obtained through calibration afterward.

\item \textbf{Achievability via simple retraining
(Theorems~\ref{thm:model-free}--\ref{thm:parametric}).}
Combining Theorem~\ref{thm:meta} with standard retraining yields
equilibrium-valid coverage when $L<1$. For discretizable classes, a
model-free variant achieves
$\widetilde{O}(1/((1{-}L)^{3/2}\sqrt{n}))$ under only $L<1$. Under
standard regularity, a parametric variant with tail-averaging after a short stabilization
period matches the $n$- and $L$-dependence of the informative minimax
rate when $\beta/\gamma = O(1)$.

\end{enumerate}

\paragraph{The price of performativity.}
The factor $1/(1{-}L)^2$ in the sample complexity is the \emph{price of
performativity}.  It quantifies how endogenous distribution shift
inflates the data required to achieve a fixed excess miscoverage level.
Relative to the non-performative setting, the required sample size grows
by a factor of order $1/(1{-}L)^2$ up to logarithmic terms.  For
example, $L=0.5$ corresponds to a $4\times$ increase, while $L=0.9$
corresponds to a $100\times$ increase.  As $L \to 1$, this cost
diverges, matching the impossibility at $L \ge 1$.

\paragraph{Paper outline.}
Section~\ref{sec:setup} formalizes equilibrium-valid coverage, the $L$-sensitivity framework, and the transfer principle.
Section~\ref{sec:main-results} establishes the fundamental limits: impossibility when $L \ge 1$, minimax lower bounds when $L < 1$, and the price of performativity interpretation.
Section~\ref{sec:algorithm} gives retraining-based procedures and their upper bounds.
Section~\ref{sec:experiments} evaluates our framework empirically, beginning with a semi-synthetic study that demonstrates the value of equilibrium-aware calibration, followed by a set of controlled experiments that support the predicted sample complexity scaling, the $L\epsilon$ transfer effect, and the phase transition at $L=1$. Section~\ref{sec:related} positions the results relative to related work on performative prediction and conformal prediction under shift.
Section~\ref{sec:discussion} discusses practical implications and extensions.

\section{Framework and the Transfer Principle}
\label{sec:setup}

We first recall the performative prediction framework and then
introduce the objects used in our analysis.

\subsection{Performative Prediction}
\label{sec:perf-setup}

Let $\cX$ be a feature space, $\cY$ a label space, and $\cF$ a class
of predictors $f: \cX \to \cY$.  Following
\citet{perdomo2020performative}, a \emph{distribution map}
$D: \cF \to \Delta(\cX \times \cY)$ assigns to each predictor $f$ a
joint distribution $D(f)$ over $(X,Y)$.  We write $D_f$ as shorthand.

\begin{definition}[Performative stability {\citep{perdomo2020performative}}]
A predictor $\fstar \in \cF$ is \emph{performatively stable} if
\[
  \fstar \in \arg\min_{f \in \cF}\,
  \E_{(X,Y)\sim D(\fstar)}\!\bigl[\ell(f(X),Y)\bigr].
\]
\end{definition}

The key parameter in our analysis is $L$-sensitivity, which captures how strongly the distribution changes when the predictor changes.

\begin{definition}[$L$-sensitivity]\label{def:sensitivity}
$D(\cdot)$ is $L$-sensitive if for all $f, g \in \cF$:
\[
\TV\bigl(D(f),\, D(g)\bigr) \;\leq\; L \cdot d_{\cF}(f, g),
\]
where $\TV$ denotes total variation distance and $d_\cF$ is a metric
on $\cF$.
\end{definition}

TV-sensitivity is the right assumption for our transfer result because
total variation satisfies
$\TV(P,Q)=\sup_A |P(A)-Q(A)|$, so it controls differences in event
probabilities uniformly over all measurable events. The coverage event
$A=\{(x,y): s(x,y;f)\le q\}$ is one such event. Therefore, to
bound the coverage gap between $D(\hat f)$ and $D(\fstar)$, it suffices
to control the total variation distance between these two distributions.

\paragraph{Standing assumptions.}
Throughout this section and Section~\ref{sec:algorithm}, we assume
the distribution map $D(\cdot)$ is $L$-sensitive
(Definition~\ref{def:sensitivity}) and that an equilibrium
$\fstar$ exists (i.e., a performatively stable predictor).
We also assume that the nonconformity score $s$ is measurable and bounded. Any
additional conditions are stated with each theorem.

\subsection{Equilibrium-Valid Coverage}
\label{sec:our-formulation}

A nonconformity score $s: \cX \times \cY \times \cF \to \R$ measures
how unusual $(x,y)$ is under predictor $f$.  The prediction set at
threshold $q$ is
$C_f(x) = \{y \in \cY : s(x,y;f) \leq q\}$.  Our theory is agnostic
to the choice of score: all results hold for any measurable $s$.

The guarantee in~\eqref{eq:equil-coverage} requires validity under the
equilibrium distribution $D(\fstar)$.  The next two definitions make
the feasibility question precise.

\begin{definition}[Nontrivial Prediction Set]\label{def:nontrivial}
A prediction set mapping $C: \cX \to 2^\cY$ is \emph{nontrivial} if
$C(x) \neq \cY$ for at least one $x \in \cX$ with positive
probability under the induced distribution.
\end{definition}

\begin{definition}[Excess Miscoverage]\label{def:excess-miscoverage}
Given a target miscoverage level $\alpha \in (0,1)$, the \emph{excess miscoverage}
of a prediction set $C_f$ at distribution $P$ is
\[
\mathrm{ExMisc}(C_f, P)
= \max\{0,\, \Pr_{P}(Y \notin C_f(X)) - \alpha\}.
\]
\end{definition}

\begin{definition}[Marginal Excess Miscoverage]\label{def:marginal-excess}
For a calibration-dependent procedure producing $C_S$ from a sample $S$,
the \emph{marginal excess miscoverage} at $P$ is
\[
\overline{\mathrm{ExMisc}}(C, P)
= \max\!\Bigl\{0,\;
\E_S\!\bigl[\Pr_{P}(Y\notin C_S(X))\bigr]-\alpha\Bigr\}.
\]
\end{definition}

\noindent
The nontriviality condition excludes the vacuous solution
$C(x)=\cY$, which achieves perfect coverage but conveys no
information.  Definition~\ref{def:nontrivial} is essentially the weakest
non-vacuity requirement: it only demands that the prediction set be
strictly smaller than $\cY$ with positive probability.  This suffices to
state the impossibility result (Theorem~\ref{thm:impossibility}), and
that result continues to hold under any stronger informativeness
criterion.  For the minimax \emph{rate} (Theorem~\ref{thm:minimax}),
however, a purely qualitative condition is too weak: on a non-singleton
$\cX$ an algorithm can output $\cY$ except on a vanishingly small region,
remain nontrivial in the sense of Definition~\ref{def:nontrivial}, and
still drive excess miscoverage to zero.  We therefore use the following
quantitative version.

\begin{definition}[$\tau$-Informativeness]\label{def:informative}
Fix $\tau\in(0,1]$.  A prediction-set procedure producing $C_S$ from a
calibration sample $S$ is \emph{$\tau$-informative} at $P$ if
\[
\mathrm{Inf}(C, P)
\;\coloneqq\;
\E_S\!\Bigl[\Pr_{X\sim P_X}\!\bigl(C_S(X)\neq\cY\bigr)\Bigr]
\;\ge\;\tau,
\]
where $P_X$ denotes the marginal distribution of $X$ under $P$.
\end{definition}

\noindent
Setting $\tau\to 0^+$ corresponds to Definition~\ref{def:nontrivial}.  Excess
miscoverage quantifies the shortfall from the target coverage level
$1-\alpha$.

\paragraph{Interaction model.}
The learning process proceeds over multiple rounds.  At round $t$, the
algorithm selects a predictor $f_t$ and observes $n_t$ independent
samples from the induced distribution $D(f_t)$.  The total sample size
is $n = \sum_t n_t$.

\paragraph{Minimax formulation.}
We evaluate performance in a minimax sense.  An algorithm interacts
with an unknown distribution map $D$ that is $L$-sensitive and admits a
stable fixed point $\fstar$.  The adversary selects both $D$ and a nonconformity score
$s : \cX \times \cY \times \cF \to [0,1]$.
The predictor class $\cF$, metric $d_\cF$, and sensitivity level $L$
are fixed and known.

Without an informativeness constraint the minimax problem is vacuous:
the constant procedure $C(x)=\cY$ has zero excess miscoverage under every
$D$ and $s$, so the unconstrained infimum is $0$.  We therefore optimize
only over procedures that are genuinely informative at equilibrium.  For
a fixed $\tau\in(0,1]$, let
\[
\mathcal A_\tau
=
\bigl\{
\text{algorithms : }
\mathrm{Inf}(C_{\hat f}, D(\fstar))\ge\tau
\text{ for every admissible } (D,s)
\bigr\}
\]
be the class of $\tau$-informative algorithms
(Definition~\ref{def:informative}).  The \emph{$\tau$-informative minimax
marginal excess miscoverage} is
\[
\mathfrak{M}_\tau(n, L, \cF)
=
\inf_{\mathrm{alg}\,\in\,\mathcal A_\tau}\, \sup_{D,\, s}
\overline{\mathrm{ExMisc}}(C_{\hat{f}}, D(\fstar)).
\]

\noindent
In words, the learner seeks an informative algorithm that minimizes
worst-case equilibrium marginal excess miscoverage over all $L$-sensitive
distribution maps and scoring rules.  The minimax value depends on $n$
and $\mathcal F$, which determine statistical difficulty, and on $L$,
which determines performative difficulty.

This objective underlies our lower bound
(Theorem~\ref{thm:minimax}) and matching upper bounds
(Theorems~\ref{thm:model-free}--\ref{thm:parametric}).

\subsection{Coverage Transfer Theorem}
\label{sec:meta}

The following lemma turns distributional mismatch into a coverage penalty.

\begin{lemma}[Coverage Decomposition]\label{lem:decomp}
Fix predictors $h,g \in \cF$. Let $S \sim D(g)^m$ be a held-out calibration
sample, and let $C_{h,S}$ be the split conformal prediction set obtained
by calibrating the scores $s(\cdot;h)$ on $S$. Then
\[
  \Pr_{S,\,(X,Y)\sim D(h)}\!\bigl(Y \in C_{h,S}(X)\bigr)
  \;\geq\; 1 - \alpha - \TV\!\bigl(D(h), D(g)\bigr).
\]
\end{lemma}
\noindent
(Proof in Appendix~\ref{app:upper}.) Split conformal calibration on $S$
guarantees $1-\alpha$ coverage for a test point drawn from $D(g)$. If the
test point is instead drawn from $D(h)$, the coverage probability changes
by at most the total variation distance between the two distributions.

The independence condition is standard in split conformal prediction:
conditional on $h$, the calibration scores on $S$ and the test score
must be exchangeable under the same distribution for the $1-\alpha$
coverage guarantee to apply. Without this independence, the split
conformal guarantee can fail because the calibration sample $S$ may
depend on $h$. In our algorithms, we enforce this by drawing $S$
\emph{after} the final predictor is fixed.

\begin{theorem}[Coverage Transfer Theorem]\label{thm:meta}
Assume the distribution map $D(\cdot)$ is $L$-sensitive with respect to
$d_\cF$, with $L<1$.
Let $\hat f$ satisfy $d_\cF(\hat f,\fstar)\le \epsilon$.
Construct $C_{\hat f}$ by split conformal calibration using a held-out
calibration sample $S \sim D(\hat f)^m$ drawn after $\hat f$ is fixed.
Then
\[
\Pr_{S,\,(X,Y)\sim D(\fstar)}\!\bigl(Y\in C_{\hat f,S}(X)\bigr)
\;\ge\; 1-\alpha - L\epsilon.
\]
Equivalently,
\[
\max\!\left\{0,\;
\Pr_{S,\,(X,Y)\sim D(\fstar)}\!\bigl(Y\notin C_{\hat f,S}(X)\bigr)-\alpha
\right\}
\;\le\; L\epsilon.
\]
\end{theorem}

\begin{proof}[Proof sketch]
\text{Step 1 (validity under $D(\hat f)$).}
Because calibration uses held-out i.i.d.\ data from $D(\hat f)$ drawn
after $\hat f$ is fixed, split conformal guarantees
\[
\Pr_{S,\,(X,Y)\sim D(\hat f)}\!\bigl(Y\in C_{\hat f,S}(X)\bigr)
\ge 1-\alpha.
\]

\text{Step 2 (transfer to $D(\fstar)$).}
For each realized calibration sample $S$, define
\[
A_S=\{(X,Y):Y\in C_{\hat f,S}(X)\}.
\]
This is the event that the prediction set constructed from $S$
covers the test label.
By total variation,
\[
\Pr_{(X,Y)\sim D(\fstar)}(A_S)
\ge
\Pr_{(X,Y)\sim D(\hat f)}(A_S)
-
\TV\!\bigl(D(\hat f),D(\fstar)\bigr).
\]
This holds for every realization of $S$, so averaging over $S$ gives
\[
\Pr_{S,\,(X,Y)\sim D(\fstar)}\!\bigl(Y\in C_{\hat f,S}(X)\bigr)
\ge
1-\alpha-\TV\!\bigl(D(\hat f),D(\fstar)\bigr).
\]
By $L$-sensitivity,
\[
\TV(D(\hat f),D(\fstar))
\le
L\,d_\cF(\hat f,\fstar)
\le
L\epsilon.
\]
Combining the two steps proves the claim.
\end{proof}

\subsection{Example: Lending}
\label{sec:example-lending}

We illustrate $L$-sensitivity in a simple lending model.

A bank chooses a credit threshold $\theta \in [0,1]$: applicants with score $X \ge \theta$ receive a loan. Receiving a loan changes future repayment behavior: access to capital provides liquidity that may decrease the probability of default, or it may increase the borrower's total debt burden and increase the probability of default. Suppose
\[
\Pr(Y=1 \mid X=x)=g(x,\theta), \quad g(x,\theta) = g_0(x) + \lambda \mathbf{1}[x \ge \theta],
\]
where $\lambda$ represents the net causal effect of receiving the loan.

The distribution of scores $X$ is fixed.
Thus, changing $\theta$ affects the data distribution only through applicants whose loan decision flips.
If $\theta$ changes to $\theta'$, the only affected applicants are those with
\[
X \in [\min(\theta,\theta'),\, \max(\theta,\theta')].
\]

The total variation distance between $D(\theta)$ and $D(\theta')$ is therefore proportional to the probability mass of this interval:
\[
\TV(D(\theta), D(\theta'))
\le
\lambda \Pr\big(X \in [\min(\theta,\theta'),\, \max(\theta,\theta')]\big).
\]

If the score distribution has density bounded by $B$, then
\[
\Pr(X \in [\theta',\theta]) \le B|\theta-\theta'|,
\]
so
\[
\TV(D(\theta), D(\theta'))
\le
\lambda B |\theta-\theta'|.
\]

Thus the distribution map is $L$-sensitive with $L=\lambda B$.

This example makes the phase transition intuitive.  When the treatment
effect $\lambda$ is large or many applicants are concentrated near the
threshold (large $B$), the environment reacts strongly to small changes
in $\theta$, corresponding to large $L$. When $\lambda B < 1$,
performative effects are moderate, and equilibrium guarantees remain
feasible.

\subsection{Strategic Classification}
\label{app:strategic-example}

Strategic classification~\citep{hardt2016strategic,dong2018revealedpreferences,milli2019socialcost,miller2020causal,levanon2021practical} provides a
natural instance of performative prediction.  A population of agents
has features $X \in \cX$ and labels $Y \in \{0,1\}$.  A classifier
$f$ is deployed.  Each agent with features $x$ can modify their
features to $x' = \Delta(x; f)$ at cost $c(x, x')$.  The
post-manipulation distribution is $D(f)$.

\paragraph{Sensitivity from manipulation costs.}
If agents' best-response functions are Lipschitz in~$f$, then
$D(\cdot)$ is $L$-sensitive with
$L = \sup_{x} \mathrm{Lip}(\Delta(x; \cdot))$.  Agents with high
manipulation costs (inelastic best-response) contribute small~$L$;
agents with low costs (elastic best-response) contribute large~$L$.

\paragraph{Interpretation.}
$L \geq 1$ means agents can shift features at least as fast as the
classifier changes. Our impossibility result (Theorem~\ref{thm:impossibility}) establishes that nontrivial prediction sets are impossible in this regime. $L < 1$ means manipulation is bounded, and
$1/(1{-}L)^2$ is the ``price of performativity'' for sample complexity. At $L = 0.7$, a practitioner needs roughly $11\times$ more data than without strategic behavior to achieve the same target coverage.

\paragraph{Worked example.}
Let $\cX = \R^d$, $\cY = \{0,1\}$, and
$f_\theta(x) = \mathbf{1}[\theta^\top x \geq 0]$.  Agents trade off the
benefit of a positive classification against the cost of manipulation,
solving the penalized best response
\[
  \Delta(x;\theta) = \arg\min_{x'}\; -\,\theta^\top x' + c\|x'-x\|^2 .
\]
Setting the gradient to zero gives $x' = x + \theta/(2c)$.
The best-response Lipschitz constant in $\theta$ therefore scales as
$1/c$: cheap manipulation (small~$c$) gives large~$L$; expensive
manipulation (large~$c$) gives small~$L$.

When $c$ is large enough that $L < 1$, the conditions for our positive
results hold, so equilibrium-valid conformal prediction becomes
possible. The resulting prediction sets may be wider than naive conformal sets at
the same sample size because they must account for the additional
uncertainty induced by performative feedback. In contrast, naive
procedures that ignore manipulation need not maintain coverage at
equilibrium.

This strategic-classification perspective motivates our
semi-synthetic experiment in
Section~\ref{sec:semi_synthetic}. There we vary the feedback strength
in a lending-related dataset to illustrate the role of $L$ in
equilibrium-valid coverage.

\section{Main Results: A Phase Transition at \texorpdfstring{$L = 1$}{L = 1}}
\label{sec:main-results}

Our main results show a sharp phase transition at $L=1$. When $L \ge 1$,
equilibrium-valid coverage is impossible. When $L < 1$, it is
achievable, and the required sample size is inflated by a factor of
$\widetilde{\Theta}(1/(1{-}L)^2)$ relative to the non-performative
setting. This section states these limits and explains their
interpretation. Matching algorithms appear in
Section~\ref{sec:algorithm}.

\subsection{Impossibility When \texorpdfstring{$L \geq 1$}{L >= 1}}
\label{sec:impossibility}

When $L \geq 1$, the feedback is too strong for uniform nontrivial
coverage guarantees.  The theorem below makes this precise.

\begin{theorem}[Impossibility for \texorpdfstring{$L \ge 1$}{L >= 1}]
\label{thm:impossibility}
Let $L \ge 1$, $\alpha \in (0,1)$, and $n \in \mathbb{N}$.
There is no algorithm that, using a total sample budget $n$, can uniformly
guarantee both of the following over all $L$-sensitive performative
environments admitting a performatively stable predictor $\fstar$:

\begin{enumerate}[label=(\alph*),itemsep=1pt]
\item \emph{Equilibrium-valid coverage:}
\[
\Pr_{(X,Y)\sim D(\fstar)}\!\bigl(Y \in C_{\hat f}(X)\bigr) \ge 1-\alpha;
\]
\item \emph{Nontriviality:}
\[
\Pr_{X \sim D_X(\fstar)}\!\bigl(C_{\hat f}(X) \neq \cY\bigr) > 0,
\]
where $D_X(\fstar)$ denotes the distribution of inputs $X$ under $D(\fstar)$.
\end{enumerate}

Equivalently, for every such algorithm, there exists an
$L$-sensitive environment for which either equilibrium-valid
coverage fails or the output prediction set is trivial.
\end{theorem}

\begin{proof}[Proof sketch]
It suffices to prove the claim for $L=1$, because every
$1$-sensitive map is automatically $L$-sensitive for every $L \ge 1$.

Consider a one-dimensional predictor class
$\cF=\Theta=[0,1]$, a singleton feature space $\cX=\{x_0\}$,
binary labels $\cY=\{0,1\}$, and squared loss
$\ell(\theta,y)=(\theta-y)^2$.
Define two $1$-sensitive distribution maps:
\begin{align*}
D_+(\theta):\quad 
Y &\sim \mathrm{Bern}\!\left(
\min\bigl(\theta+\epsilon,\;1-\epsilon\bigr)
\right), \\
D_-(\theta):\quad 
Y &\sim \mathrm{Bern}\!\left(
\max\bigl(\theta-\epsilon,\;\epsilon\bigr)
\right).
\end{align*}
These maps admit performatively stable predictors,
namely $\theta_+^\star = 1-\epsilon$ and $\theta_-^\star = \epsilon$.
Hence the corresponding equilibrium distributions are
$\mathrm{Bern}(1-\epsilon)$ and $\mathrm{Bern}(\epsilon)$.

At any deployment $\theta$, the Bernoulli parameters under
$D_+(\theta)$ and $D_-(\theta)$ differ by at most $2\epsilon$, and
both remain in $[\epsilon,1-\epsilon]$.  Therefore the per-sample KL
divergence is $O(\epsilon)$ uniformly over $\theta$.  Choosing
$\epsilon = c/n$ keeps the total KL over $n$ adaptive
samples bounded by a constant, so no algorithm can reliably
distinguish the two environments.

However, because $\cX$ is a singleton and $\cY=\{0,1\}$, the only
nontrivial prediction sets are $\{0\}$ and $\{1\}$.
If $\epsilon < \min(\alpha,1-\alpha)$, then under
$D_+(\theta_+^\star)=\mathrm{Bern}(1-\epsilon)$ the only nontrivial
set achieving coverage at least $1-\alpha$ is $\{1\}$, while under
$D_-(\theta_-^\star)=\mathrm{Bern}(\epsilon)$ the only nontrivial
set achieving coverage at least $1-\alpha$ is $\{0\}$.
Thus no algorithm can simultaneously guarantee both equilibrium-valid
coverage and nontriviality uniformly over the two environments.
By Le Cam's two-point testing lemma~\citep{lecam1986asymptotic}, the bounded KL divergence implies that no algorithm can distinguish $D_+$ from $D_-$ with high enough probability to reliably output the distinct valid sets required by each environment.

The full proof, including the extension to general spaces, is in
Appendix~\ref{app:impossibility}. That extension shows that this
binary construction can be embedded into an arbitrary feature and
label space.

\end{proof}

The implication is direct. When $L \geq 1$, collecting more data cannot yield uniform guarantees, 
because each deployment can induce a distribution shift as large as the 
information gained. The remedy is to reduce~$L$ below one, rather than 
increasing the sample size.

\subsection{Minimax Rate When \texorpdfstring{$L < 1$}{L < 1}}
\label{sec:rates}

When $L < 1$, the distribution map is a contraction and
equilibrium-valid coverage becomes feasible. In this regime,
we establish a minimax lower bound that characterizes the
optimal sample complexity.

\begin{theorem}[Minimax Lower Bound for \texorpdfstring{$L < 1$}{L < 1}]\label{thm:minimax}
Let $D(\cdot)$ be $L$-sensitive with $L < 1$, and fix any
$\tau\in(0,1]$.  For $n$ total samples, the $\tau$-informative minimax
marginal excess miscoverage satisfies
\[
  \mathfrak{M}_\tau(n, L, \cF)
  \;\geq\;
  \widetilde{\Omega}\!\left(\frac{1}{(1{-}L)\sqrt{n}}\right).
\]
Equivalently, for every algorithm in $\mathcal A_\tau$ there exists an
$L$-sensitive map such that
\[
  \overline{\mathrm{ExMisc}}(C_{\hat{f}}, D(\fstar))
  \;\geq\;
  \widetilde{\Omega}\!\left(\frac{1}{(1{-}L)\sqrt{n}}\right).
\]
The restriction to $\mathcal A_\tau$ rules out vacuous and near-vacuous
prediction sets; the rate itself is independent of $\tau$, in that any
fixed $\tau\in(0,1]$ yields the same order.
\end{theorem}

\noindent
Matching upper bounds appear in Section~\ref{sec:algorithm}.  We sketch
the lower bound here; the full proof is in Appendix~\ref{app:lower}.

\paragraph{Lower bound construction.}
We use Le~Cam's two-point method. Consider two $L$-sensitive maps on
$\cX=\{x_0\}$ and $\cY=\{0,1\}$:
\[
D_i(\theta):\quad Y \sim \mathrm{Bern}(a + i\delta + L\theta),
\qquad i\in\{0,1\},
\]
where $a$ is bounded away from $0$ and $1$.
Because the risk minimizer under squared loss is the mean, the performative equilibria satisfy
\[
\theta_i^\star = \frac{a+i\delta}{1-L}.
\]
Thus the two environments differ by only $\delta$ at any fixed
deployment $\theta$, but their equilibria are separated by
\[
|\theta_1^\star-\theta_0^\star| = \frac{\delta}{1-L}.
\]
This is the key amplification effect: a local perturbation of size
$\delta$ becomes an equilibrium gap larger by a factor $1/(1-L)$.

At the same time, the two environments are statistically hard to distinguish:
\[
p_1(\theta)-p_0(\theta)=\delta
\quad\text{for every }\theta,
\]
so the per-sample KL divergence is $O(\delta^2)$ uniformly over the
deployed predictor. Adaptivity therefore does not help: no matter
which predictor the algorithm deploys, each sample provides only
$O(\delta^2)$ information about which environment is active. Over $n$
samples, the total KL is $O(n\delta^2)$.

Choosing $\delta=c/\sqrt{n}$ keeps the total KL bounded by a constant,
so Le~Cam's lemma implies that the two environments cannot be reliably
distinguished. By centering the target coverage threshold perfectly between the two equilibrium parameters, a prediction set valid for one environment is strictly invalid for the other. Therefore, failing to distinguish the environments forces the algorithm to output an invalid set, incurring an excess miscoverage of at least
\[
\Omega\!\left(\frac{\delta}{1-L}\right)
=
\Omega\!\left(\frac{1}{(1-L)\sqrt{n}}\right).
\]
The full proof is given in Appendix~\ref{app:lower}.

\subsection{The Price of Performativity}
\label{sec:price}

Theorems~\ref{thm:impossibility} and~\ref{thm:minimax} together give a
complete picture as a function of~$L$:

\begin{center}
\begin{tabular}{@{}ll@{}}
\toprule
\textbf{Form} & \textbf{Expression} \\
\midrule
Excess miscoverage given $n$ samples
& $\widetilde{\Theta}\!\left(\dfrac{1}{(1{-}L)\sqrt{n}}\right)$\rlap{$^{\dagger}$} \\[8pt]
Sample complexity for target excess $\epsilon$
& $\widetilde{\Theta}\!\left(\dfrac{1}{(1{-}L)^2\,\epsilon^2}\right)$\rlap{$^{\dagger}$} \\
\bottomrule
\end{tabular}
\end{center}

\noindent
{\small ${}^\dagger$Upper bound achieved by the parametric algorithm
under standard regularity ($\beta/\gamma = O(1)$); the model-free
algorithm achieves $(1{-}L)^{-3/2}$ in place of $(1{-}L)^{-1}$.}

\noindent
The factor $1/(1{-}L)^2$ is the \emph{price of performativity}: how the
sample size must grow to maintain a fixed accuracy under endogenous
distribution shift.  If excess miscoverage scales as
$\widetilde{O}\!\big(1/((1{-}L)\sqrt{n})\big)$, then requiring it to be at most
$\epsilon$ implies
\[
\widetilde{O}\!\Big(\frac{1}{(1{-}L)\sqrt{n}}\Big)\le \epsilon
\quad\Longrightarrow\quad
n \ge \widetilde{\Omega}\!\Big(\frac{1}{(1{-}L)^2\epsilon^2}\Big).
\]
Thus the cost is constant at $L=0$ and grows rapidly as $L\to 1$
(e.g., $L=0.5$ requires about $4\times$ more samples, while $L=0.9$
requires about $100\times$).  This matches the divergence suggested by
the impossibility at $L=1$.

\noindent
The intuition is that performative feedback magnifies statistical error.
A local estimation error of size $\epsilon$ does not remain of size
$\epsilon$ at equilibrium. Instead, because the predictor changes the
distribution and the distribution feeds back into the predictor,
controlling equilibrium error requires paying an amplification factor on
the order of $1/(1-L)$. Therefore, to guarantee final excess miscoverage
at most $\Delta$, the underlying statistical error must be correspondingly
smaller, on the order of $\Delta(1-L)$. Since achieving estimation
accuracy $\epsilon$ requires a sample size of order $1/\epsilon^2$, this
leads to sample complexity scaling as $1/(\Delta^2(1-L)^2)$. Thus, at a
high level, the $(1-L)^{-2}$ dependence reflects the interaction between
feedback amplification and the standard quadratic cost of statistical
precision.

\section{Algorithms}
\label{sec:algorithm}

Section~\ref{sec:meta} reduces equilibrium-valid coverage to a
constructive problem: produce a predictor close to performative
equilibrium, then transfer that predictor error into a coverage
guarantee.  This section gives two retraining-based algorithms that
follow this template.  The model-free method gives general
achievability under minimal assumptions, while the parametric method
sharpens the rate under stronger regularity.

\subsection{Model-Free Algorithm}
\label{sec:model-free}

We begin with a finite-class model-free setting, meaning that we do not
assume a parametric or differentiable form for the performative map. The main structural
assumption is $L<1$, and we do not impose smoothness or strong
convexity on the performative map. In this setting, repeated retraining
combined with a final split conformal calibration step yields the
following upper bound.

\begin{algorithm}[ht]
\caption{Model-Free Performative Conformal Prediction (MF-PCP)}
\label{alg:mfpcp}
\begin{algorithmic}[1]
\REQUIRE Finite predictor class $\cF$ with $|\cF| = N$,
  non-conformity score function $s$, loss function $\ell$,
  sample budget $n$, rounds $T$, target miscoverage level $\alpha$,
  initial predictor $f_0 \in \cF$
\STATE Set $m = \lfloor n / (T+1) \rfloor$
\FOR{$t = 1, \ldots, T$}
  \STATE Deploy $f_{t-1}$ and collect $m$ i.i.d.\ samples $S_t$
    from $D(f_{t-1})$
  \STATE $f_t \gets \arg\min_{f \in \cF} \widehat{R}(f; S_t)$
    \hfill\COMMENT{ERM over $\cF$ using loss $\ell$}
\ENDFOR
\STATE Deploy $f_T$ and collect $m$ held-out samples $S_{T+1}$
  from $D(f_T)$
\STATE $\hat{q} \gets \lceil(1-\alpha)(|S_{T+1}|+1)\rceil$-th
  smallest score on $S_{T+1}$ using $s(x,y;f_T)$
\RETURN $f_T$, $C_{f_T}(x) = \{y : s(x,y;f_T) \le \hat{q}\}$
\end{algorithmic}
\end{algorithm}

For a finite class with bounded loss, empirical risk minimization (ERM)
based on $m$ samples achieves excess risk
$O(\sqrt{\log N / m})$~\citep{shalevshwartz2014understanding}.
Applied at round $t$, this controls the excess risk of $f_t$ under the
current distribution $D(f_{t-1})$. To convert such risk control into
control of the predictor itself, we assume an error-bound condition:
there exists $\mu>0$ such that, for the relevant round-wise population
risk,
\[
  d_\cF(f,f^\dagger_{t-1})
  \le
  \frac{1}{\mu}\Bigl(R(f;D(f_{t-1})) - R(f^\dagger_{t-1};D(f_{t-1}))\Bigr),
\]
where $f^\dagger_{t-1}$ denotes the population risk minimizer under
$D(f_{t-1})$.
Under this assumption, ERM yields a round-wise statistical error
\[
  d_\cF(f_t,f^\dagger_{t-1})
  =
  O\!\left(\frac{1}{\mu}\sqrt{\frac{\log N}{m}}\right),
\]
which can then be combined with the $L$-sensitivity contraction to
bound $d_\cF(f_t,\fstar)$.

\begin{theorem}[Model-Free Upper Bound]
\label{thm:model-free}
Under $L$-sensitivity with $L<1$, bounded loss $\ell \in [0,1]$,
$|\cF| = N$, and error-bound parameter $\mu > 0$,
Algorithm~\ref{alg:mfpcp} achieves approximate equilibrium-valid
coverage. More precisely, for $0<L<1$, if
\[
  T = \left\lceil \frac{\log n}{\log(1/L)} \right\rceil,
\]
then
\[
  \overline{\mathrm{ExMisc}}(C_{f_T,S}, D(\fstar))
  \;\le\;
  \widetilde{O}\!\left(\frac{1}{(1-L)^{3/2}\sqrt{n}}\right).
\]
For $L=0$, the same conclusion holds with $T=1$.
\end{theorem}

\begin{proof}[Proof sketch]
Let $m=\lfloor n/(T+1)\rfloor$ be the per-round sample size. By the ERM
bound and the error-bound condition, each retraining step incurs a
statistical error of order
\[
  \widetilde{O}_{N,\mu}\!\left(\frac{1}{\sqrt{m}}\right)
\]
in $d_\cF$. Repeated retraining therefore gives
\begin{align*}
d_\cF(f_T,\fstar)
&\le
L^T d_\cF(f_0,\fstar)
+
\sum_{k=0}^{T-1} L^k\,
\widetilde{O}_{N,\mu}\!\left(\frac{1}{\sqrt{m}}\right) \\
&\le
L^T d_\cF(f_0,\fstar)
+
\frac{1}{1-L}\,
\widetilde{O}_{N,\mu}\!\left(\frac{1}{\sqrt{m}}\right).
\end{align*}
The first term is the transient error from the initialization, while
the second is the accumulated statistical error across rounds.

With
\[
  T = \left\lceil \frac{\log n}{\log(1/L)} \right\rceil,
\]
we have $L^T \le 1/n$, so the transient term is negligible. Moreover,
since $\log(1/L) \ge 1-L$ for $0<L<1$, the number of rounds satisfies
\[
  T = O\!\left(\frac{\log n}{1-L}\right),
\]
and hence
\[
  m = \Omega\!\left(\frac{n(1-L)}{\log n}\right).
\]
Substituting this into the accumulated error term yields
\[
  \frac{1}{1-L}\,
  \widetilde{O}_{N,\mu}\!\left(\frac{1}{\sqrt{m}}\right)
  \;=\;
  \widetilde{O}_{N,\mu}\!\left(\frac{1}{(1-L)^{3/2}\sqrt{n}}\right).
\]
Thus,
\[
  d_\cF(f_T,\fstar)
  \;=\;
  \widetilde{O}_{N,\mu}\!\left(\frac{1}{(1-L)^{3/2}\sqrt{n}}\right).
\]

The retraining phase bounds the distance between $f_T$ and $\fstar$.
The independent calibration sample then ensures conformal validity
under $D(f_T)$. Theorem~\ref{thm:meta} translates this distance bound
into the stated excess-miscoverage guarantee under $D(\fstar)$.
\end{proof}

\begin{remark}
The model-free rate has an extra $(1-L)^{-1/2}$ factor relative to
the minimax lower bound, arising from splitting the sample budget
across $O(\log n/(1-L))$ rounds.  In the parametric setting,
tail-averaging after a short stabilization period closes this gap
(Theorem~\ref{thm:parametric})
by exploiting the AR(1) autocorrelation structure of the iterates.
Whether the information-theoretic rate is achievable over finite
classes without smoothness assumptions is an open question.
\end{remark}

\begin{remark}[Beyond Finite Classes]
Although Theorem~\ref{thm:model-free} is stated for a finite class
$\cF$, the same type of guarantee extends to many infinite classes.
One can approximate $\cF$ by a finite subset that is dense enough so
that the approximation error is of smaller order than the statistical
error, and then run the same procedure on that subset.  This modifies
the bound only by additional logarithmic factors that reflect the size
of the approximation. 
\end{remark}

\subsection{Parametric Algorithm}
\label{sec:parametric}

We now turn to continuous parametric classes, where predictors are
parameterized by vectors in a convex subset of $\mathbb{R}^d$. In this
setting, discretization is less attractive because it can introduce
substantial dimension-dependent error. Instead, under stronger
regularity assumptions, we analyze repeated empirical risk
minimization directly.

Under these assumptions, the population retraining map is contractive,
so repeated retraining drives the predictor toward equilibrium.
Finite samples introduce noise at each round.  Averaging later
iterates reduces this noise, but averaging too early retains too much
dependence on the initial predictor.  We therefore use a short
stabilization period and average only the later iterates, in the spirit
of Polyak--Ruppert averaging~\citep{ruppert1988efficient,polyak1992acceleration}.  

\begin{algorithm}[ht]
\caption{Parametric Performative Conformal Prediction (P-PCP)}
\label{alg:pcp}
\begin{algorithmic}[1]
\REQUIRE Initial predictor $f_0$, total sample budget $n$,
stabilization rounds $B$, averaging rounds $K$, target miscoverage level $\alpha$
\STATE Set $m = \lfloor n/(B+K+1) \rfloor$
\FOR{$t = 1, \ldots, B+K$}
  \STATE Deploy $f_{t-1}$ and collect $m$ samples $S_t \sim D(f_{t-1})$
  \STATE Retrain by ERM:
  \[
  f_t \gets \arg\min_{f\in\cF} \widehat R(f;S_t)
  \]
\ENDFOR
\STATE Form the tail-averaged predictor
\[
\bar f_{B,K} \gets \frac{1}{K}\sum_{t=B+1}^{B+K} f_t
\]
\STATE Deploy $\bar f_{B,K}$ and collect $m$ fresh samples
$S_{B+K+1} \sim D(\bar f_{B,K})$ for calibration
\STATE Let $\hat q$ be the
$\lceil (1-\alpha)(m+1)\rceil$-th smallest score in
$\{s(x,y;\bar f_{B,K}):(x,y)\in S_{B+K+1}\}$
\RETURN $\bar f_{B,K}$ and
\[
C_{\bar f_{B,K}}(x)=\{y:s(x,y;\bar f_{B,K})\le \hat q\}
\]
\end{algorithmic}
\end{algorithm}

\begin{theorem}[Parametric Upper Bound]
\label{thm:parametric}
Assume $\cF \subseteq \mathbb{R}^d$ is convex, with
$d_\cF(\cdot,\cdot) = \|\cdot - \cdot\|_2$.
Assume $D(\cdot)$ is $L$-sensitive.
Assume also that for every distribution $D$, the population risk $R_D$
is $\beta$-smooth and $\gamma$-strongly convex over $\cF$.
Let
\[
\rho \coloneqq \frac{L\beta}{\gamma},
\qquad \rho < 1.
\]
Run Algorithm~\ref{alg:pcp} with stabilization length and averaging
window satisfying
\[
B = O\!\left(\frac{\log n}{1-\rho}\right),
\qquad
K = O\!\left(\frac{\log n}{1-\rho}\right).
\]
Then
\[
\overline{\mathrm{ExMisc}}(C_{\bar f_{B,K}},D(\fstar))
\;\le\;
\widetilde O\!\left(\frac{1}{(1-\rho)\sqrt n}\right).
\]
When $\rho \asymp L$, for example in the well-conditioned regime
$\beta/\gamma \approx 1$, this matches the $n$- and $L$-dependence of the
informative minimax lower bound (Theorem~\ref{thm:minimax}) up to
logarithmic factors under the stated regularity conditions.
\end{theorem}

\noindent
\emph{Proof sketch.}
Under smoothness and strong convexity, the population minimizer varies
Lipschitzly with the underlying distribution; combined with
$L$-sensitivity, this implies that the population retraining map is
contractive with factor $\rho=L\beta/\gamma$.

Thus repeated retraining moves the predictor toward equilibrium, and
the effect of the initial predictor shrinks by a factor $\rho$ at each
round.
We choose the stabilization length $B$ so that this initialization term
is no larger than the final statistical error.
After this, the remaining $K$ iterates fluctuate around equilibrium due
to sampling noise, and averaging them reduces that noise.
This yields
\[
d_\cF(\bar f_{B,K},\fstar)
=
\widetilde O\!\left(\frac{1}{(1-\rho)\sqrt n}\right).
\]
The transfer theorem (Theorem~\ref{thm:meta}) then converts this
predictor error into the stated excess miscoverage bound.
Full details are given in Appendix~\ref{app:upper}.

\begin{remark}[Upper bounds are $\tau$-informative]\label{rem:informative-upper}
The lower bound of Theorem~\ref{thm:minimax} is stated over the
constrained class $\mathcal A_\tau$, so a matching upper bound is only
meaningful for procedures that are themselves $\tau$-informative.  Both
algorithms output a split-conformal set at the finite empirical quantile
$\hat q$, the $\lceil(1-\alpha)(m+1)\rceil$-th smallest calibration
score.  Whenever the equilibrium score distribution is non-degenerate
above its $(1-\alpha)$-quantile---a mild condition beyond the standing
boundedness assumption---the quantile $\hat q$ stays bounded away from
$\sup s$, so the returned set excludes at least one label with
probability bounded below; that is,
$\mathrm{Inf}(C_{\hat f},D(\fstar))\ge\tau$ for some $\tau>0$.  Both
algorithms then lie in $\mathcal A_\tau$, so the upper and lower bounds
compare within the same algorithm class.
\end{remark}

\subsection{Comparing the Two Algorithms}
\label{sec:comparison}

Table~\ref{tab:comparison} summarizes the tradeoffs between the two approaches. The model-free method requires only bounded loss but incurs an additional $(1-L)^{-1/2}$ penalty. The parametric method assumes smoothness and strong convexity, allowing retraining to contract toward equilibrium. By tail-averaging to reduce sampling noise, the parametric approach achieves a sharper rate, removing the extra $(1-L)^{-1/2}$ factor (up to logarithmic terms) when $\beta/\gamma \approx 1$.

\begin{table}[h]
\centering
\caption{Comparison of the two algorithms}
\label{tab:comparison}
\setlength{\tabcolsep}{3pt}
\renewcommand{\arraystretch}{1.1}
\begin{tabularx}{\columnwidth}{@{}l>{\RaggedRight\arraybackslash}X>{\RaggedRight\arraybackslash}X@{}}
\toprule
& \textbf{Model-free} (Alg.~\ref{alg:mfpcp})
& \textbf{Parametric} (Alg.~\ref{alg:pcp}) \\
\midrule
Contraction
& $L<1$
& $\rho = L\beta/\gamma < 1$ \\
Output
& Last iterate
& Average after stabilization \\
Rate
& $\widetilde{O}\!\bigl(\frac{1}{(1-L)^{3/2}\sqrt{n}}\bigr)$
& $\widetilde{O}\!\bigl(\frac{1}{(1-\rho)\sqrt{n}}\bigr)$ \\
Loss assumptions
& Bounded
& Smooth, strongly convex \\
Relation to lower bound
& Suboptimal by a factor $(1-L)^{-1/2}$, up to logs
& Matches lower bound up to logs when $\rho \asymp L$ \\
\bottomrule
\end{tabularx}
\end{table}

This comparison also clarifies the mechanics of the transfer theorem. The minimax lower bound depends exclusively on the performative sensitivity $L$ (Theorem~\ref{thm:minimax}), capturing the statistical difficulty induced by feedback. By contrast, the parametric upper bound depends on $\rho=L\beta/\gamma$, coupling the feedback strength $L$ with the optimization geometry ($\beta$ and $\gamma$) of the loss.
\section{Experiments}
\label{sec:experiments}

This section serves two purposes.  We first show practical relevance in
a semi-synthetic credit-scoring task, where performative feedback
causes a standard conformal baseline to fail and where equilibrium-aware
retraining restores coverage.  We then turn to controlled synthetic
models that isolate the dependence predicted by the theory on $L$, $n$,
and predictor error.

\subsection{Semi-Synthetic Case Study: Strategic Regression in Credit Scoring}
\label{sec:semi_synthetic}

We begin with a semi-synthetic strategic regression task based on the
Kaggle ``Give Me Some Credit'' dataset~\citep{kaggle2011givemesomecredit}.  The goal is to test the method
under realistic covariates and model misspecification, rather than only
in stylized synthetic settings.

\begin{figure}[tbp]
    \centering

    \begin{subfigure}[t]{0.49\linewidth}
        \centering
        \includegraphics[width=\linewidth]{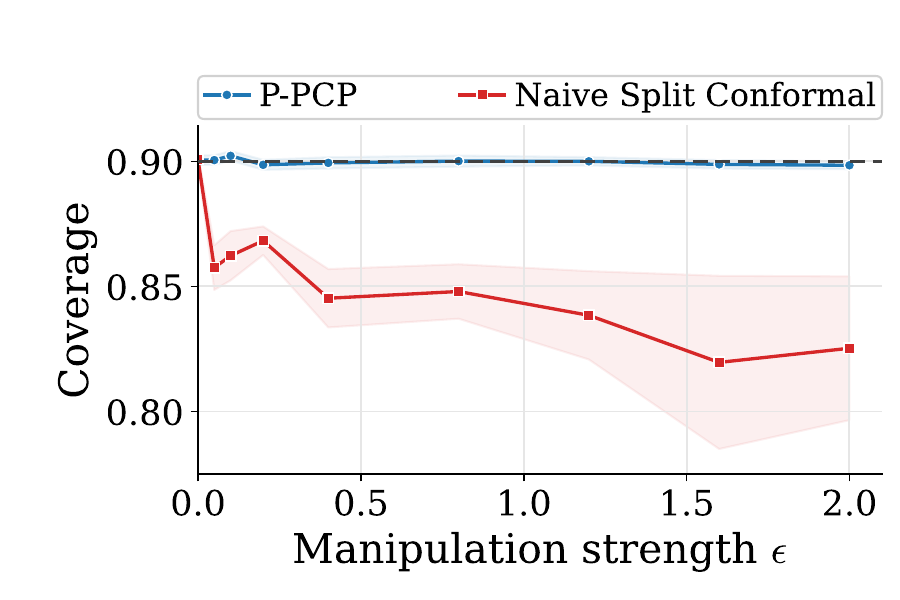}
        \vspace{-1.5em}
        \caption{Equilibrium coverage versus $\epsilon$}
        \label{fig:semi_syn_cov}
    \end{subfigure}\hfill
    \begin{subfigure}[t]{0.49\linewidth}
        \centering
        \includegraphics[width=\linewidth]{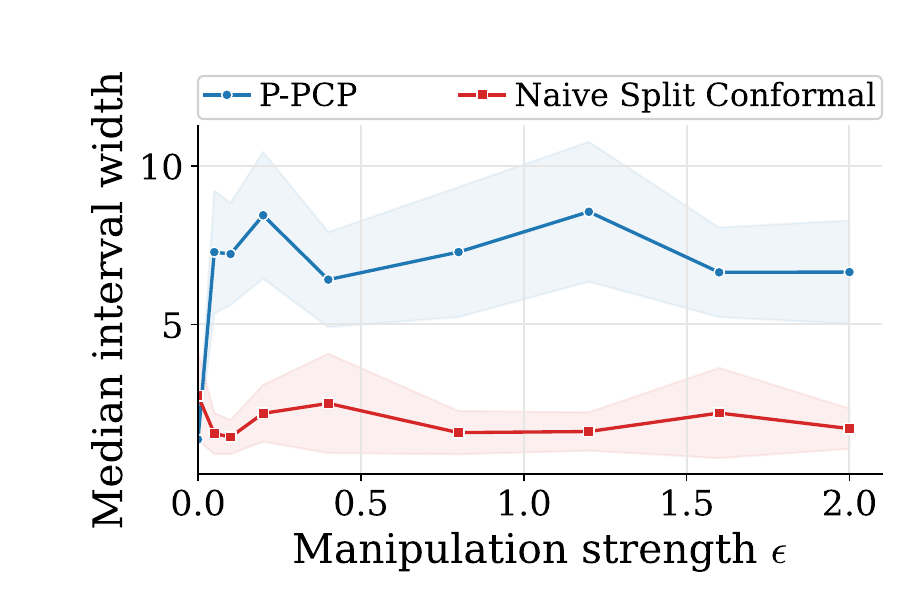}
        \vspace{-1.5em}
        \caption{Prediction interval width versus $\epsilon$}
        \label{fig:semi_syn_width}
    \end{subfigure}

    \caption{Semi-synthetic strategic regression experiment.
    (a) Equilibrium coverage versus manipulation strength $\epsilon$.
    The naive split conformal baseline undercovers as performativity increases,
    whereas P-PCP remains closer to the nominal $90\%$ target.
    (b) Median prediction interval width versus $\epsilon$ for the two procedures.}
    \label{fig:semi_synthetic_results}
\end{figure}

\paragraph{Results.}
Figure~\ref{fig:semi_synthetic_results} presents the results of the semi-synthetic experiment.
As the manipulation strength $\epsilon$ increases, the Naive Split Conformal baseline
substantially undercovers, dropping well below the nominal $90\%$ level.
Because this baseline calibrates on historical, unmanipulated data, it fails to account
for the endogenous distribution shift induced by deployment.
In contrast, P-PCP remains close to the target coverage across the full range of $\epsilon$.

Figure~\ref{fig:semi_synthetic_results}\subref{fig:semi_syn_width}
shows the corresponding interval widths. In this run, the naive baseline
stays nearly flat as $\epsilon$ increases, while P-PCP adapts by using
materially wider intervals once performative effects become non-negligible.
At high manipulation strengths, P-PCP produces markedly wider intervals
than the naive baseline.
In this semi-synthetic example, equilibrium-aware retraining with
fresh calibration preserves nominal coverage much better than
historical calibration alone.

\subsection{Synthetic Illustrations of the Theory}

We now turn to controlled synthetic models.  Their role is different:
they isolate the dependence on $L$, $n$, and predictor error so that
the scaling laws from the theory can be checked directly.  Unless
stated otherwise, we use the Gaussian location model
$Y \sim \cN(a + L\theta, 1)$ with $a = 0.3$, score
$s(y;\theta)=|y-\theta|$, and target miscoverage level $\alpha=0.1$.  Each
configuration uses 2000 independent replications and 50{,}000 test
samples drawn from $D(\fstar)$.  Full simulation settings are deferred
to Appendix~\ref{app:experiments} in the supplementary material.

\begin{figure}[tbp]
    \centering

    \begin{subfigure}[t]{0.95\linewidth}
        \centering
        \includegraphics[width=\linewidth]{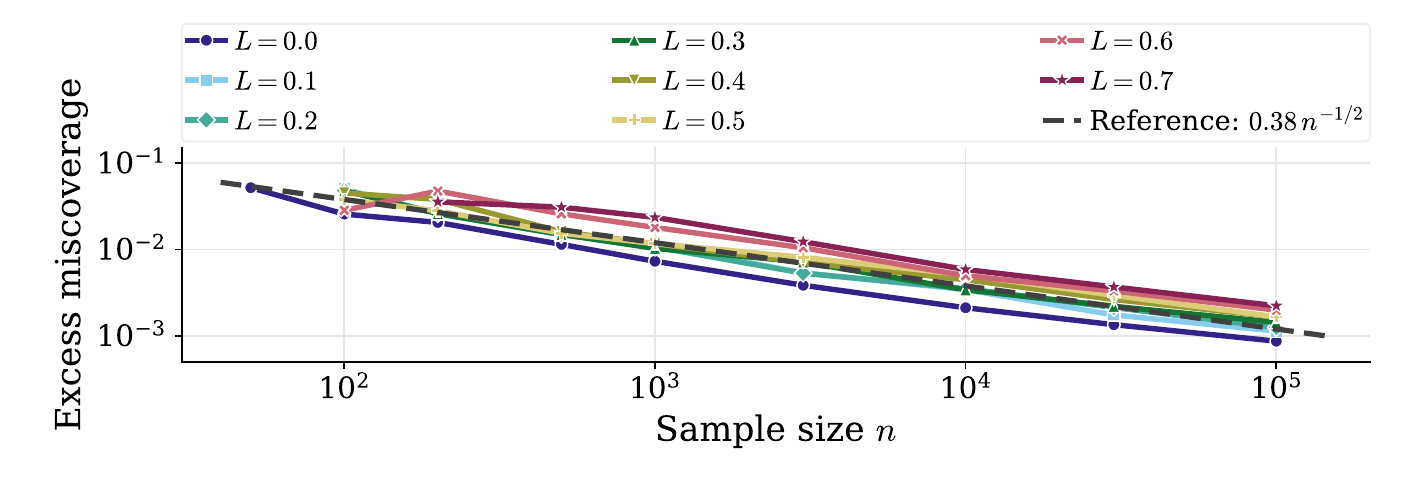}
        \vspace{-1.5em} %
        \caption{Excess miscoverage versus sample size $n$}
        \label{fig:rate}
    \end{subfigure}

    \vspace{0.75em} %

    \begin{subfigure}[t]{0.95\linewidth}
        \centering
        \includegraphics[width=\linewidth]{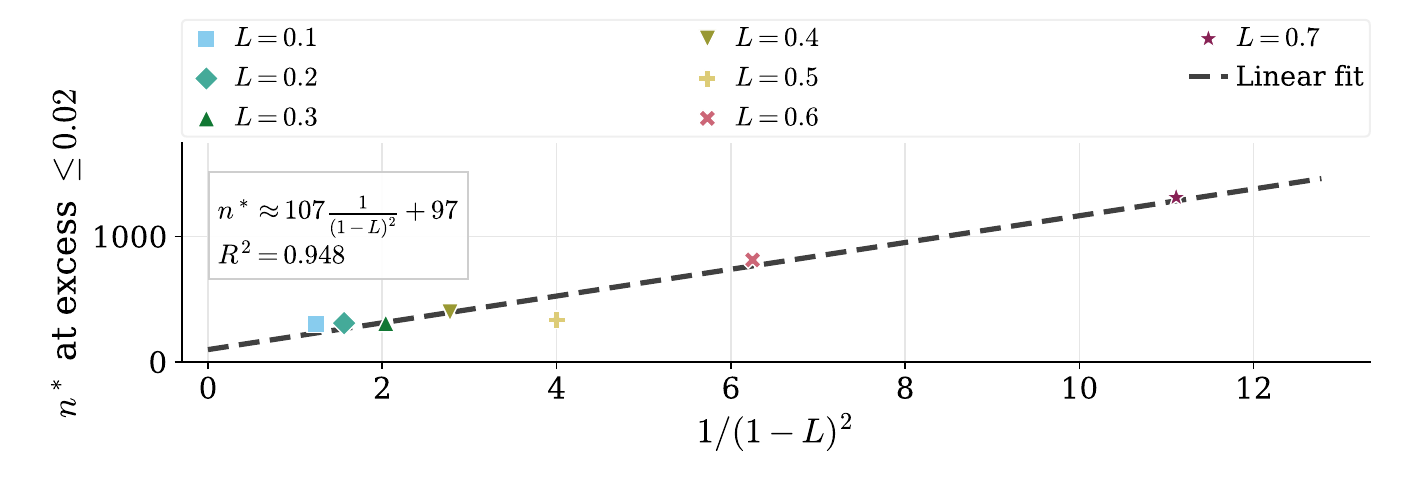}
        \vspace{-1.5em}
        \caption{Threshold sample size $n^\star(L)$}
        \label{fig:nstar}
    \end{subfigure}
    \vspace{0.75em}
    \begin{subfigure}[t]{0.95\linewidth}
        \centering
        \includegraphics[width=\linewidth]{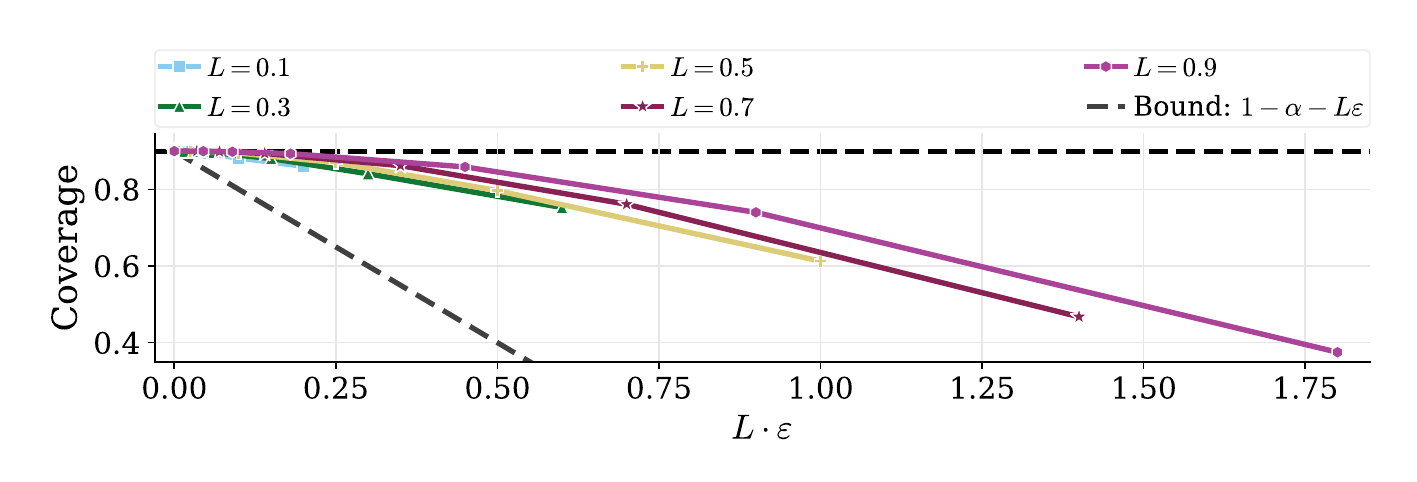}
        \vspace{-1.5em}
        \caption{Coverage versus $L\varepsilon$}
        \label{fig:collapse}
    \end{subfigure}

    \caption{Synthetic validation of the theory.
    (a) Excess miscoverage versus sample size $n$ (log--log) for
    $L\in\{0,\ldots,0.7\}$; the dashed line shows the $n^{-1/2}$
    reference slope.
    (b) Threshold sample size $n^\star(L)$, defined as the smallest
    $n$ for which excess miscoverage is at most $0.02$, plotted against
    $1/(1-L)^2$; the dashed line is a linear fit.
    (c) Coverage versus $L\varepsilon$ for five values of $L$.
    The dashed line is the first-order lower bound
    $1-\alpha-L\varepsilon$.}
    \label{fig:synthetic_all}
\end{figure}

\paragraph{Rate and price of performativity
(Simulation~1; Theorem~\ref{thm:minimax}).}
We use Algorithm~\ref{alg:pcp} as a testbed to examine the scaling
predicted by Theorem~\ref{thm:minimax}.  For each fixed $L$, the theorem
predicts excess miscoverage of order
$\widetilde{\Theta}(1/((1-L)\sqrt{n}))$, so the dependence on sample
size is $n^{-1/2}$, while the dependence on performative sensitivity
appears through the prefactor $(1-L)^{-1}$.  We vary
$L \in \{0,0.1,\ldots,0.7\}$ and
$n \in \{50,\ldots,100{,}000\}$.
Figure~\ref{fig:synthetic_all}\subref{fig:rate} plots excess
miscoverage against $n$ on a log--log scale. The curves are approximately parallel to the
$n^{-1/2}$ reference slope (dashed), which matches the predicted decay
in sample size.  At the same time, larger values of $L$ shift the
curves upward, consistent with the predicted deterioration as $L$
approaches $1$.

Figure~\ref{fig:synthetic_all}\subref{fig:nstar} shows the same
relationship in sample-complexity form. Since Theorem~\ref{thm:minimax} predicts excess
miscoverage of order $\widetilde{\Theta}(1/((1-L)\sqrt{n}))$,
the sample size needed to reach a fixed target excess should scale as
$\widetilde{\Theta}(1/(1-L)^2)$.  To test this prediction, we define
$n^\star(L)$ as the smallest sample size for which the excess
miscoverage falls below $0.02$.  The horizontal axis plots the
predicted complexity term $1/(1-L)^2$, and the vertical axis plots the
corresponding empirical threshold $n^\star(L)$.  The points lie close
to a straight line, with fitted relation
$n^\star \approx 107/(1-L)^2 + 97$
and $R^2 = 0.948$, in line with the predicted
$\widetilde{\Theta}(1/(1-L)^2)$ scaling.  The nonzero intercept
reflects the baseline cost already present in the non-performative
case.

\paragraph{Transfer theorem
(Simulation~2; Theorem~\ref{thm:meta}).}
We next test the main prediction of the transfer theorem.  For a
deployed predictor at distance
\[
\varepsilon = |\theta_{\mathrm{deployed}}-\theta^\star|
\]
from equilibrium, the theorem implies that the leading-order coverage
loss should depend on the single quantity $L\varepsilon$.  We fix
$m=2000$ calibration samples and vary $\varepsilon$ across five values
of $L$.

Figure~\ref{fig:collapse} plots coverage against $L\varepsilon$.  Near
the origin, the curves align closely: for $L\varepsilon \le 0.3$, their
spread is at most $0.03$.  This is the regime where the first-order
prediction is expected to be most accurate.  Over this range, all
curves remain above the lower bound $1-\alpha-L\varepsilon$ (dashed).
The figure therefore supports the prediction that, for small predictor
error, coverage degradation is governed primarily by the product
$L\varepsilon$, rather than by $L$ and $\varepsilon$ separately.

\begin{figure}[tbp]
    \centering

    \begin{subfigure}[t]{0.95\linewidth}
        \centering
        \includegraphics[width=\linewidth]{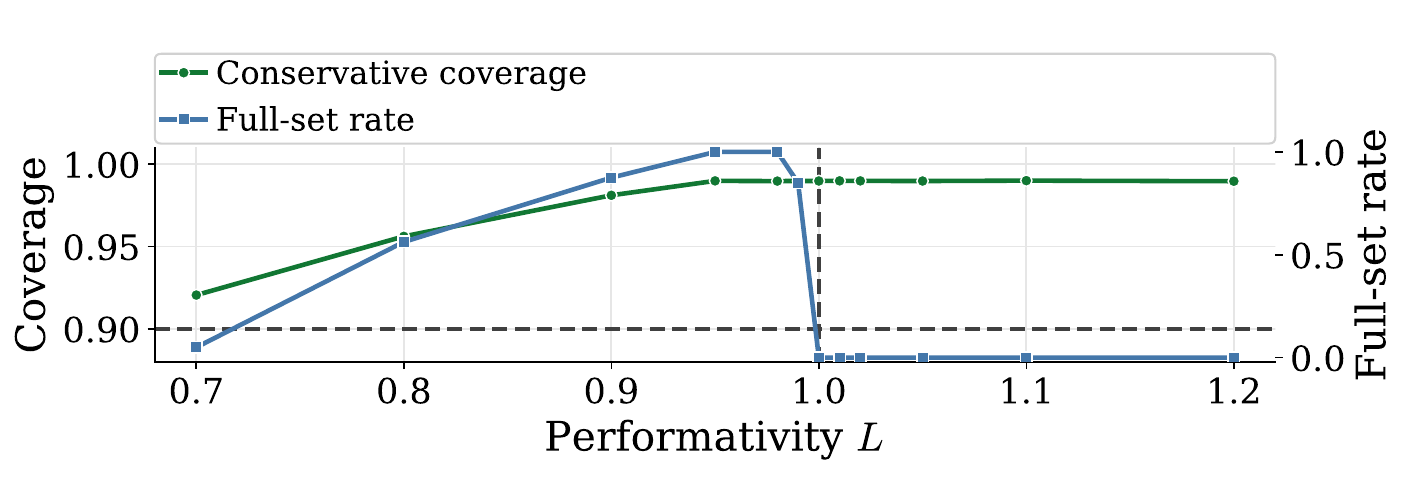}
        \vspace{-1.5em}
        \caption{Bounded Bernoulli sweep}
        \label{fig:phase_transition_unrestricted}
    \end{subfigure}

    \caption{Phase transition at $L=1$.
    Bounded Bernoulli sweep based on
    $g_{L,\pm}(\theta)=\operatorname{clip}(c+L(\theta-c)\pm\varepsilon,\,
    [\varepsilon,1-\varepsilon])$ with $c=0.05$ and
    $\varepsilon=0.01$. The green curve is conservative equilibrium
    coverage (minimum over the $+$ and $-$ branches); the blue curve is
    the unrestricted full-set rate (maximum over the two branches).
    The supercritical regime is therefore well-defined for all plotted
    $L$, with explicit equilibria throughout. The dashed vertical line
    marks the threshold $L=1$.}
    \label{fig:phase_transition}
\end{figure}

\paragraph{Phase transition at $L=1$
(Simulation~3; Theorem~\ref{thm:impossibility}).}
We next illustrate the feasibility boundary predicted by
Theorem~\ref{thm:impossibility} using bounded Bernoulli maps whose
equilibria remain explicit on both sides of the threshold.
Figure~\ref{fig:phase_transition} reports the results.

Figure~\ref{fig:phase_transition} sweeps
$L \in \{0.7,0.8,0.9,0.95,0.98,0.99,1.0,1.01,1.02,1.05,1.1,1.2\}$ for
the clipped family
\[
g_{L,\pm}(\theta)=\operatorname{clip}(c+L(\theta-c)\pm\varepsilon,\,
[\varepsilon,1-\varepsilon]),
\qquad c=0.05,\ \varepsilon=0.01,
\]
and reports conservative equilibrium coverage (the minimum over the two
branches) together with the unrestricted full-set rate (the maximum over
the two branches). As $L$ approaches $1$ from below, coverage remains
near nominal only while the unrestricted procedure increasingly returns
the full label set. Unlike the previous Gaussian illustration, the
supercritical regime here is still a genuine equilibrium experiment,
because the clipped Bernoulli maps are well-defined and have explicit
fixed points for every plotted $L$.
The figure is therefore a theorem-compatible empirical illustration of
the feasibility boundary in a setting that remains well-defined across
the full sweep, including $L \ge 1$. As before, it is an empirical
illustration, not an empirical proof of impossibility.

\section{Related Work}
\label{sec:related}

\paragraph{Performative prediction.}
The performative prediction framework of \citet{perdomo2020performative}
shows that when $L<1$, repeated retraining converges to a stable
equilibrium, and that convergence slows as $L \to 1$. Subsequent work
studies stochastic retraining dynamics~\citep{mendler2020stochastic},
regret under performativity~\citep{jagadeesan2022regret},
state-dependent feedback~\citep{brown2022stateful}, and performative
power~\citep{hardt2022performative}. This entire line of work asks how
quickly a system reaches equilibrium and how to minimize performative
risk. We ask a different question: once the system has stabilized, how
reliably can we quantify prediction uncertainty? The connection is that
the same parameter $L$ that governs convergence also determines whether
valid uncertainty guarantees are achievable.

\paragraph{Conformal prediction under distribution shift.}
Several papers extend conformal prediction to settings where the test
distribution differs from the calibration distribution.
\citet{tibshirani2019covariateshift} handle covariate shift via
importance-weighted scores. \citet{barber2023conformal} derive coverage
bounds for general non-exchangeable data.
\citet{gibbs2021adaptive,gibbs2024online} give procedures that maintain coverage under arbitrary temporal drift.
\citet{oliveira2024nonexchangeable} study split conformal guarantees
when data is not exchangeable. None of these papers model the
performative case, where deployment itself alters the future
data-generating process.

The closest papers on the conformal side are
\citet{fannjiang2022conformal} and \citet{binny2025who}.
\citeauthor{fannjiang2022conformal} study a biomolecular design setting
in which deploying a model changes which candidate molecules are proposed
in the next round. However, their analysis handles only covariate shift:
the marginal distribution over inputs changes, but the conditional
distribution of outcomes given inputs does not. In a performative
system, predictions affect behavior, which changes outcomes, so the full
joint distribution shifts. Moreover, their analysis does not address the
equilibrium that emerges after repeated retraining.
\citeauthor{binny2025who} apply a similar iterative conformal idea to
multi-agent motion planning, where deploying an uncertainty-aware
controller changes how surrounding agents behave. Their focus is on
safe trajectory prediction rather than performative equilibrium in the
statistical learning sense, and they do not analyze sample complexity.
Our work is the first to study conformal validity at the equilibrium
distribution $D(\fstar)$ reached by repeated retraining.

\paragraph{Uncertainty quantification under performativity.}
Three recent papers study uncertainty quantification in performative
settings, and we differ from each in a distinct way.

\citet{li2025prc} introduce Performative Risk Control (PRC), which
calibrates a scalar threshold so that a chosen risk measure stays
controlled throughout the retraining process. Their guarantee is that
$\E[\ell(z;\lambda)] \le \alpha$ at each step, which is not the same as
constructing a prediction set that covers the true outcome with
probability $1-\alpha$: scalar risk control does not in general imply
set-valued coverage. They also use a different sensitivity notion,
measuring how much the expected loss changes when the threshold changes,
rather than how much the data distribution itself changes. As a result, their analysis does not extend to coverage guarantees for
prediction sets, does not yield a transfer theorem connecting proximity
to equilibrium with coverage validity, and provides no impossibility
results or minimax lower bounds.

\citet{csillag2024strategic} construct prediction sets that are valid
against the worst-case strategic response of the population. Because
they do not assume anything about which equilibrium will emerge, their
sets must be wide enough to cover all possible outcomes under any
plausible response, and may be unnecessarily wide when the actual
equilibrium is less extreme than the worst case. We instead target the
specific equilibrium $D(\fstar)$ reached by repeated retraining, which
gives tighter sets when $L < 1$. The two approaches suit different
situations: theirs when the response is unknown and robustness is the
priority; ours when the system stabilizes and one can calibrate at the
stable distribution. Theorem~\ref{thm:impossibility} delineates where
each applies: when $L < 1$, targeting the equilibrium is possible and
the sample complexity cost is $\widetilde{\Theta}(1/(1-L)^2)$
times the non-performative baseline; when $L \ge 1$, no procedure can
simultaneously achieve coverage and nontrivial set sizes.

\citet{perdomo2025multical} shows that performative multicalibration --- simultaneously matching predicted probabilities to observed outcome frequencies across a large collection of subgroups --- is always achievable under performative feedback, without any Lipschitz condition on $D(\cdot)$. This is a non-trivial algorithmic result, since multicalibration is demanding even in the static setting. Since multicalibration is a form of validity guarantee, one might read this as evidence that useful uncertainty quantification is always achievable under performativity, even when $L \ge 1$, which would contradict our impossibility result. It does not, because the two results target different objects. Multicalibration is a property of predicted probabilities: it asks that the model's outputs average correctly within subgroups. Coverage is a property of prediction sets: it asks that the set $C(x)$ excludes some labels while still containing the true outcome with probability $1-\alpha$. A predictor can be perfectly multicalibrated and yet produce no useful prediction sets, because calibration of probabilities does not constrain which labels a set includes or excludes. When $L \ge 1$, it is this exclusion requirement that becomes impossible to satisfy at equilibrium.

No prior work establishes the feasibility threshold for
equilibrium-valid coverage, derives the minimax sample complexity below
that threshold, proves a transfer theorem that converts proximity to
equilibrium directly into a coverage guarantee, or provides algorithms
that match the optimal rate. These are the four contributions of the
present paper.

\section{Discussion and Conclusion}
\label{sec:discussion}
Our results establish a sharp phase transition at $L=1$ for conformal
prediction under performative feedback: nontrivial equilibrium-valid
coverage is impossible when $L \ge 1$, and achievable with sample
complexity inflated by a factor of $\widetilde{\Theta}(1/(1-L)^2)$
when $L < 1$. We now discuss practical implications, limitations,
and extensions.

\paragraph{Practical implications.}
Our results suggest a workflow for uncertainty quantification
under performativity.  Before deployment in an environment where
predictions affect outcomes, estimate the sensitivity $L$, for example
from pilot deployments.  If $L \geq 1$, then equilibrium-valid
prediction sets are not uniformly achievable: the main obstacle is an
unstable feedback loop, not a lack of data.  In that regime, the
relevant intervention is to reduce $L$ by limiting how strongly the
predictor changes decisions.  If $L < 1$, then the non-performative sample requirement should be inflated by a factor
of at least $1/(1-L)^2$ (and up to $1/(1-L)^3$ without smoothness and
strong convexity) to budget for valid equilibrium coverage. The transfer theorem further
shows that any improvement in estimating or reaching equilibrium
translates directly into better equilibrium-valid coverage guarantees.

\paragraph{Limitations.}
Our analysis has several limitations.  First, the coverage guarantees
hold on average at equilibrium across the whole population, rather than
separately for each feature value or subgroup.  Second, the sharp
parametric upper bound requires smoothness and strong convexity, while
the more general model-free result pays an additional
$(1-L)^{-1/2}$ factor.  Third, the experiments are primarily synthetic or semi-synthetic studies meant to test the predicted scaling laws and illustrate equilibrium-aware calibration, rather than full application-specific evaluations. Finally, the theory represents performative feedback
through an $L$-sensitive distribution map, which summarizes the
strength of the feedback in a single parameter.  Extending the results
to more local notions of sensitivity is a natural next step.

\paragraph{Extensions.}
Our minimax characterization also points to natural next directions.
The model-free algorithm (Theorem~\ref{thm:model-free}) shows that
$L < 1$ is both necessary and sufficient for finite or compact
predictor classes, while the parametric algorithm achieves the sharper
rate for continuous classes under the standard regularity condition
$\rho = L\beta/\gamma < 1$.  A central open question is whether the
information-theoretic rate is achievable in polynomial time under only
$L < 1$, without strong convexity.  The transfer theorem
(Theorem~\ref{thm:meta}) is useful precisely because it decouples this
optimization question from calibration: any future method with a
convergence guarantee to equilibrium immediately inherits a coverage
guarantee through the theorem.

\paragraph{Reproducibility.}
The synthetic and semi-synthetic experiments rely only on the models and
parameter grids described in Section~\ref{sec:experiments} and
Appendix~\ref{app:experiments}.  The semi-synthetic study uses the
publicly available ``Give Me Some Credit'' dataset
\citep{kaggle2011givemesomecredit}.  Code and simulation scripts to
reproduce all figures will be made publicly available upon acceptance.

\bibliography{refs}

\onecolumn
\appendix
\section{Proof of Theorem~\ref{thm:impossibility}}
\label{app:impossibility}

It suffices to prove the theorem for $L=1$, because if
\[
\TV(D(f),D(g)) \le d_\cF(f,g)
\quad\text{for all } f,g\in\cF,
\]
then automatically
\[
\TV(D(f),D(g)) \le L\, d_\cF(f,g)
\quad\text{for all } f,g\in\cF
\]
for every $L \ge 1$. Thus any hard $1$-sensitive instance is also a
hard instance for every $L \ge 1$.

\subsection{A Binary Hard Instance for Theorem~\ref{thm:impossibility}}

\begin{proof}
Take $\cX = \{x_0\}$, $\cY = \{0,1\}$, $\Theta=[0,1]$, and squared loss
$\ell(\theta,y)=(\theta-y)^2$. Fix $\epsilon>0$, to be chosen as a
function of $n$. Define two distribution maps:
\begin{align*}
D_+(\theta) &: Y \sim \mathrm{Bern}\!\bigl(g_+(\theta)\bigr),
\qquad
g_+(\theta) = \min(\theta+\epsilon,\; 1-\epsilon),\\
D_-(\theta) &: Y \sim \mathrm{Bern}\!\bigl(g_-(\theta)\bigr),
\qquad
g_-(\theta) = \max(\theta-\epsilon,\; \epsilon).
\end{align*}

\textbf{Step 1: Both maps are $1$-sensitive.}
We first show that $g_+$ and $g_-$ are $1$-Lipschitz.
Write
\[
g_+(\theta)=c(\theta+\epsilon),
\qquad
c(u)=\min(u,1-\epsilon).
\]
The map $\theta \mapsto \theta+\epsilon$ is $1$-Lipschitz, since
\[
|(\theta+\epsilon)-(\theta'+\epsilon)|=|\theta-\theta'|.
\]
The clipping map $c$ is also $1$-Lipschitz, because truncation at a
constant cannot increase distances:
\[
|c(u)-c(v)| \le |u-v|
\qquad\text{for all }u,v\in\mathbb{R}.
\]
Hence $g_+$, as a composition of two $1$-Lipschitz maps, is
$1$-Lipschitz. The same argument applies to
\[
g_-(\theta)=\max(\theta-\epsilon,\epsilon),
\]
so $g_-$ is also $1$-Lipschitz.

For Bernoulli distributions, total variation distance is the absolute
difference of the success probabilities. Therefore
\[
\TV(D_\pm(\theta),D_\pm(\theta'))
=
|g_\pm(\theta)-g_\pm(\theta')|
\le
|\theta-\theta'|.
\]
Hence both distribution maps are $1$-sensitive.

The clipping ensures that the Bernoulli parameter remains in $[\epsilon, 1-\epsilon]$,
which is needed both for the distribution to be well-defined and for
the KL divergence bound in Step~3.  Without clipping, the parameter
could exceed $1$ or fall below $0$; any other truncation that preserves
the $1$-Lipschitz property (for example, clipping to $[0,1]$) would
work equally well.

\textbf{Step 2: Both maps admit performatively stable predictors.}
Under squared loss, the risk minimizer under a Bernoulli distribution
with mean $p$ is $\theta=p$. Therefore a performatively stable
predictor solves
\[
\theta^\star = \E[Y \mid Y \sim D(\theta^\star)] = g(\theta^\star).
\]

For $D_+$, the fixed-point equation is
\[
\theta = \min(\theta+\epsilon,\; 1-\epsilon).
\]
Because $\epsilon > 0$, the branch $\theta=\theta+\epsilon$ is impossible. Therefore, the fixed point must lie in the clipped region, so
\[
\theta_+^\star = 1-\epsilon.
\]
At $\theta=1-\epsilon$,
\[
g_+(1-\epsilon)=\min(1,\;1-\epsilon)=1-\epsilon,
\]
so $1-\epsilon$ is a fixed point.

For $D_-$, the fixed-point equation is
\[
\theta = \max(\theta-\epsilon,\; \epsilon),
\]
and similarly the only solution is
\[
\theta_-^\star = \epsilon,
\]
since
\[
g_-(\epsilon)=\max(0,\;\epsilon)=\epsilon.
\]

Thus the equilibrium distributions are
\[
D_+(\theta_+^\star)=\mathrm{Bern}(1-\epsilon),
\qquad
D_-(\theta_-^\star)=\mathrm{Bern}(\epsilon).
\]

\textbf{Step 3: The two environments are hard to distinguish.}
For every deployment $\theta \in [0,1]$, the two Bernoulli parameters
differ by at most $2\epsilon$. On the interior interval
$[2\epsilon,\,1-2\epsilon]$ we have
\[
g_+(\theta)=\theta+\epsilon,
\qquad
g_-(\theta)=\theta-\epsilon,
\]
and therefore
\[
|g_+(\theta)-g_-(\theta)| = 2\epsilon.
\]
If $\theta \notin [2\epsilon,\,1-2\epsilon]$, then either
$\theta-\epsilon < \epsilon$ or $\theta+\epsilon > 1-\epsilon$,
so one of the two expressions is truncated by the clipping operator.
Such truncation cannot increase the distance between them. Hence,
for all $\theta \in [0,1]$,
\[
|g_+(\theta)-g_-(\theta)| \le 2\epsilon.
\]

Moreover, both parameters lie in $[\epsilon,1-\epsilon]$. Hence
\[
\KL\!\bigl(D_+(\theta)\,\|\,D_-(\theta)\bigr)
\le
\frac{(g_+(\theta)-g_-(\theta))^2}{g_-(\theta)(1-g_-(\theta))}
\le
\frac{(2\epsilon)^2}{\epsilon(1-\epsilon)}
=
\frac{4\epsilon}{1-\epsilon}.
\]
Now consider any adaptive algorithm that deploys predictors
$\theta_1,\dots,\theta_T$ and collects $n_t$ samples at round $t$,
with $\sum_{t=1}^T n_t=n$. By the chain rule for KL divergence,
\[
\KL(\Pr_+^{(n)} \,\|\, \Pr_-^{(n)})
=
\sum_{t=1}^T
n_t \,
\E_+\!\left[
\KL\!\bigl(D_+(\theta_t)\,\|\,D_-(\theta_t)\bigr)
\right]
\le
\frac{4n\epsilon}{1-\epsilon}.
\]
Choosing
\[
\epsilon = \frac{c_0}{n}
\]
for some constant $c_0>0$ makes the total KL bounded
by a constant. Therefore Le Cam's lemma implies that no
algorithm can reliably distinguish whether the data came from
$D_+$ or $D_-$.

\textbf{Step 4: The nontrivial valid sets are incompatible.}
Because $\cX$ is a singleton and $\cY=\{0,1\}$, any prediction set is
equivalent to one of
\[
\emptyset,\qquad \{0\},\qquad \{1\},\qquad \cY.
\]
The empty set is never valid, and $\cY$ is trivial. So the only
nontrivial candidates are $\{0\}$ and $\{1\}$.

Their equilibrium coverages are:
\begin{center}
\begin{tabular}{@{}lcc@{}}
\toprule
& $D_+(\theta_+^\star)=\mathrm{Bern}(1-\epsilon)$
& $D_-(\theta_-^\star)=\mathrm{Bern}(\epsilon)$ \\
\midrule
$C=\{1\}$ & $1-\epsilon$ & $\epsilon$ \\
$C=\{0\}$ & $\epsilon$ & $1-\epsilon$ \\
\bottomrule
\end{tabular}
\end{center}

If
\[
\epsilon < \min(\alpha,\;1-\alpha),
\]
then under $D_+$ the only nontrivial set with coverage at least
$1-\alpha$ is $\{1\}$, while under $D_-$ the only nontrivial set with
coverage at least $1-\alpha$ is $\{0\}$.

Therefore any algorithm that uniformly guaranteed both
equilibrium-valid coverage and nontriviality would, from its output,
implicitly identify which environment is active: output $\{1\}$ under
$D_+$ and output $\{0\}$ under $D_-$. But Step 3 shows that no
algorithm can reliably distinguish the two environments. This is a
contradiction.

Hence, for $L=1$, no algorithm can uniformly guarantee both
equilibrium-valid coverage and nontriviality.

Every $1$-sensitive map is also $L$-sensitive for every $L \ge 1$.
Hence the same impossibility result holds for all $L \ge 1$.
\end{proof}

\subsection{Extension to general spaces}

\begin{lemma}[Reduction to the binary instance]
\label{lem:app-reduction-imp}
Assume $\cX \neq \emptyset$ and $|\cY| \ge 2$.
Pick any $x_0 \in \cX$ and any distinct labels $y_0,y_1 \in \cY$.
Consider the binary hard instance from the proof above, and view it as
a distribution on $(\cX,\cY)$ by placing all feature mass on $x_0$ and
all label mass on $\{y_0,y_1\}$, with binary label $0$ identified with
$y_0$ and binary label $1$ identified with $y_1$.
Then the binary construction embeds unchanged into the general spaces,
so the impossibility result for the binary case applies to the general
case as well.
\end{lemma}

\begin{proof}
The larger spaces $\cX$ and $\cY$ do not matter, because the embedded
construction uses only the single feature point $x_0$ and the two
labels $y_0,y_1$. All other feature values and labels have probability
zero.

Therefore the embedded instance is the same as the binary one,
just with relabeled outcomes. The total variation distances are the
same, the fixed-point equations are the same, and the coverage of any
prediction set depends only on whether it contains $y_0$ and/or $y_1$
at $x_0$. Thus the same contradiction proved in the binary case applies
unchanged.
\end{proof}

\section{Full Lower Bound Proof (Theorem~\ref{thm:minimax}, Lower Bound)}
\label{app:lower}

We prove the lower bound by constructing two performative environments
that are statistically hard to distinguish but for which a single
informative prediction set is valid at one equilibrium and strictly
invalid at the other.

\subsection{Construction}

Take $\cX=\{x_0\}$ and $\cY=\{0,1\}$.
Let the predictor be a scalar $\theta\in[0,1]$, and use squared loss
$\ell(\theta,y)=(\theta-y)^2$.
Fix $L<1$, and define two $L$-sensitive maps:
\begin{align*}
D_0(\theta) &: \quad Y \sim \mathrm{Bern}\bigl(a + L\theta\bigr),\\
D_1(\theta) &: \quad Y \sim \mathrm{Bern}\bigl(a + \delta + L\theta\bigr),
\end{align*}
where $a\in(0,1)$ is chosen away from the boundary so that the Bernoulli
parameters stay bounded away from $0$ and $1$, and $\delta>0$
is a small perturbation.
Assume $a+\delta+L<1$ so that both Bernoulli parameters remain in $(0,1)$.

This construction changes the performative map only by a constant shift
$\delta$, while keeping the slope equal to $L$.

For Bernoulli outcomes under squared loss, the population risk minimizer
is the Bernoulli mean.
Therefore the equilibrium points satisfy
\[
\theta_0^\star = a + L\theta_0^\star,
\qquad
\theta_1^\star = a + \delta + L\theta_1^\star,
\]
hence
\[
\theta_0^\star = \frac{a}{1-L},
\qquad
\theta_1^\star = \frac{a+\delta}{1-L}.
\]
Since under squared loss the population minimizer equals the Bernoulli
mean, the equilibrium Bernoulli parameters are
\[
p_0^\star = \frac{a}{1-L},
\qquad
p_1^\star = \frac{a+\delta}{1-L},
\]
so the equilibrium gap is
\[
|p_1^\star-p_0^\star|
=
\frac{\delta}{1-L}.
\]

This is the amplification effect used in the lower bound: a
local perturbation of size $\delta$ in the performative map produces an
equilibrium shift larger by a factor $1/(1-L)$.

\subsection{Sequential KL divergence}

An algorithm may deploy predictors adaptively.
At round $t$, it chooses $\theta_t$ based on the past and then receives
$n_t$ independent samples from either $D_0(\theta_t)$ or
$D_1(\theta_t)$, with $\sum_t n_t = n$.

For every deployment $\theta$,
\[
p_1(\theta)-p_0(\theta)=\delta,
\]
since $p_1(\theta)=a+\delta+L\theta$ and $p_0(\theta)=a+L\theta$.
Thus, the separation between the two Bernoulli models is the same at
every deployment point.

Because $a$ is chosen away from the boundary and $a+\delta+L<1$, both
$p_0(\theta)$ and $p_1(\theta)$ lie in a fixed interval
$[\eta,1-\eta] \subset (0,1)$, uniformly over $\theta$.

For Bernoulli distributions,
\[
\KL(\mathrm{Bern}(p)\|\mathrm{Bern}(q))
\le \frac{(p-q)^2}{q(1-q)}.
\]
Since $p_0(\theta),p_1(\theta)\in[\eta,1-\eta]$ uniformly over
$\theta$, we obtain
\[
\KL\!\bigl(\mathrm{Bern}(p_0(\theta))
\,\|\,\mathrm{Bern}(p_1(\theta))\bigr)
\le \frac{\delta^2}{\eta(1-\eta)}
=: C_a\delta^2.
\]
 
Now condition on the history up to round $t-1$. The deployment
$\theta_t$ is then fixed, and the $n_t$ samples collected at round $t$
are independent draws from either
$\mathrm{Bern}(p_0(\theta_t))$ or $\mathrm{Bern}(p_1(\theta_t))$.
Therefore the KL contribution of round $t$, conditional on the past, is
\[
n_t \,
\KL\!\bigl(
\mathrm{Bern}(p_0(\theta_t))
\,\|\,\mathrm{Bern}(p_1(\theta_t))
\bigr).
\]
Applying the chain rule for KL divergence and then taking expectation
over the history gives
\[
\KL(\Pr_0^{(n)} \| \Pr_1^{(n)})
=
\sum_{t=1}^{T} n_t \,
\E_0\!\left[
\KL\!\bigl(
\mathrm{Bern}(p_0(\theta_t))
\,\|\,\mathrm{Bern}(p_1(\theta_t))
\bigr)
\right]
\le
C_a \sum_{t=1}^T n_t \delta^2
=
C_a n \delta^2.
\]

Therefore adaptivity does not help: no matter how the algorithm chooses
its deployments, the total information available for distinguishing the
two environments is at most of order $n\delta^2$.

\subsection{From equilibrium separation to excess miscoverage}

We now explain why this equilibrium gap forces a lower bound on excess
miscoverage.

Because $\cX$ is a singleton and $\cY=\{0,1\}$, every prediction set is
equivalent to one of
\[
\emptyset,\qquad \{0\},\qquad \{1\},\qquad \{0,1\}.
\]
The empty set is always invalid and $\{0,1\}$ is trivial, so the only
nontrivial candidates are $\{0\}$ and $\{1\}$.  Since $\cX$ is a
singleton, $\mathrm{Inf}\in\{0,1\}$, so membership in $\mathcal A_\tau$
(for any $\tau\in(0,1]$) is exactly the exclusion of the trivial set
$\{0,1\}$.  Combined with validity---the empty set is never
valid---the output must be $\{0\}$ or $\{1\}$.  This is what rules out
the vacuous full-set procedure that would otherwise attain zero excess
miscoverage.

Under environment $i\in\{0,1\}$, the equilibrium label distribution is
$\mathrm{Bern}(p_i^\star)$.
Therefore:
\[
\Pr_i\!\bigl(Y\notin \{1\}\bigr)=1-p_i^\star,
\qquad
\Pr_i\!\bigl(Y\notin \{0\}\bigr)=p_i^\star.
\]
So the excess miscoverage of the two nontrivial sets is determined
directly by the equilibrium Bernoulli parameter $p_i^\star$.

To make the informative set $\{1\}$ straddle its validity boundary,
choose the target miscoverage level $\alpha$ and the base parameter $a$
so that the two equilibrium means lie on opposite sides of the threshold
$1-\alpha$.
Concretely, take
\[
p_0^\star = 1-\alpha-\frac{\Delta}{2},
\qquad
p_1^\star = 1-\alpha+\frac{\Delta}{2},
\]
where
\[
\Delta = \frac{\delta}{1-L}.
\]
This can be achieved by choosing
\[
a = (1-L)\!\left(1-\alpha-\frac{\Delta}{2}\right).
\]

Under $D_0$, the set $\{1\}$ has miscoverage
\[
1-p_0^\star = \alpha+\frac{\Delta}{2},
\]
so its excess miscoverage is $\Delta/2$.
Under $D_1$, the set $\{1\}$ has miscoverage
\[
1-p_1^\star = \alpha-\frac{\Delta}{2},
\]
so its excess miscoverage is $0$.
Because the two equilibrium means $p_0^\star,p_1^\star$ straddle the
threshold $1-\alpha$, the informative set $\{1\}$ is the one whose
validity flips across the two environments: it is \emph{valid} under
$D_1$ (excess $0$) and \emph{strictly invalid} under $D_0$ (excess
$\Delta/2$).  The other informative set $\{0\}$ does not help, since its
miscoverage is $p_i^\star$, which equals $1-\alpha\mp\Delta/2$ and hence
incurs excess of order $1-2\alpha$ under both environments.  Indeed,
under $D_0$ \emph{no} informative set is valid: the best informative
choice, $\{1\}$, still misses the target by $\Delta/2$.

Therefore the same informative set $\{1\}$ is valid at one equilibrium
and strictly invalid at the other.  Any $\tau$-informative algorithm that
cannot distinguish the two environments must commit to essentially the
same output in both; outputting $\{1\}$ (the only viable informative set)
then incurs excess miscoverage at least of order $\Delta$ under $D_0$,
while the trivial full set is excluded from $\mathcal A_\tau$.

\subsection{Completing the lower bound}

The two equilibria straddle the validity boundary of the informative set
$\{1\}$: under $D_0$ every informative set incurs excess miscoverage at
least $\Delta/2$ (the minimum, attained by $\{1\}$), whereas under $D_1$
the set $\{1\}$ is valid with zero excess. Therefore any $\tau$-informative
procedure that fails to distinguish the two environments incurs excess
miscoverage at least $\Delta/2$ in one of them.

Le Cam's two-point lemma now applies. Since the KL divergence between
the two induced data-generating processes is at most $C_a n\delta^2$,
we have
\[
\inf_{\phi}
\bigl(\Pr_0(\phi=1)+\Pr_1(\phi=0)\bigr)
\;\ge\;
1-\sqrt{\tfrac12\,C_a n\delta^2}.
\]
Hence, at least one of the two testing errors is at least
\[
\frac12\left(1-\sqrt{\tfrac12\,C_a n\delta^2}\right).
\]
Combining this with the loss separation $\Delta/2$, we obtain
\[
\max_{i\in\{0,1\}}
\overline{\mathrm{ExMisc}}(C_{\hat f}, D_i(\fstar))
\;\ge\;
\frac{\Delta}{4}
\left(
1-\sqrt{\tfrac12\, C_a n\delta^2}
\right).
\]

Set
\[
\delta = \frac{c'}{\sqrt n}
\]
with $c'>0$ small enough that the square-root term is bounded away from
$1$.
Then
\[
\Delta
=
\frac{\delta}{1-L}
=
\frac{c'}{(1-L)\sqrt n},
\]
and therefore
\[
\max_{i\in\{0,1\}}
\overline{\mathrm{ExMisc}}(C_{\hat f}, D_i(\fstar))
\;\ge\;
\Omega\!\left(\frac{1}{(1-L)\sqrt n}\right).
\]

This proves the lower bound in Theorem~\ref{thm:minimax}.

\section{Full Upper Bound Proof}
\label{app:upper}

\subsection{Proof of Coverage Decomposition (Lemma~\ref{lem:decomp})}

\begin{proof}
Define score CDFs:
\begin{align*}
F_h(q) &= \Pr_{(X,Y)\sim D(h)}(s(X,Y;h) \leq q), \\
G_g(q) &= \Pr_{(X,Y)\sim D(g)}(s(X,Y;h) \leq q).
\end{align*}

\textbf{Step 1: Conformal guarantee under the calibration distribution.}
Conditional on $h$, the scores $s_1, \ldots, s_m$ are i.i.d.\ from
$G_g$.  A new score $s_{m+1}$ drawn from the same distribution is
exchangeable with $s_1, \ldots, s_m$. By the exchangeability property of order statistics, the rank of
$s_{m+1}$ among $s_1,\ldots,s_{m+1}$ is uniformly distributed over
$\{1,\ldots,m+1\}$, meaning that each rank occurs with probability
$1/(m+1)$. Therefore,
\[
\Pr(s_{m+1} \leq \hat{q})
\;\geq\;
\frac{\lceil(1-\alpha)(m+1)\rceil}{m+1}
\;\geq\; 1 - \alpha.
\]
So $\E_S[G_g(\hat{q})] \geq 1 - \alpha$.

\textbf{Step 2 (Transfer to the target distribution).}
For any measurable event $A$:
$|P(A) - Q(A)| \leq \TV(P, Q)$.
Taking $A = \{(x,y) : s(x,y;h) \leq q\}$:
\[
|F_h(q) - G_g(q)| \leq \TV(D(h), D(g)) \quad \forall q.
\]
This holds for the random threshold $\hat{q}$:
$F_h(\hat{q}) \geq G_g(\hat{q}) - \TV(D(h), D(g))$.

\text{Step 3 (Combine).}
\begin{align*}
\E_S[F_h(\hat{q})]
&\geq \E_S[G_g(\hat{q})] - \TV(D(h), D(g)) \\
&\geq (1-\alpha) - \TV(D(h), D(g)).
\end{align*}
\end{proof}

\subsection{Proof of the Coverage Transfer Theorem (Theorem~\ref{thm:meta})}

\begin{proof}
Apply Lemma~\ref{lem:decomp} with $h = g = \hat{f}$: calibrate and
deploy the same predictor.  Split conformal on $m$ i.i.d.\ samples
from $D(\hat{f})$ gives coverage $\geq 1 - \alpha$ under
$D(\hat{f})$; this is exact in finite samples, with the conformal
quantile being conservative (tending toward overcoverage, not
undercoverage).

Transfer to $D(\fstar)$ via the TV bound: for any event~$A$,
\[
  |\Pr_{D(\hat{f})}(A) - \Pr_{D(\fstar)}(A)|
  \;\leq\; \TV(D(\hat{f}),\, D(\fstar))
  \;\leq\; L \cdot d_\cF(\hat{f},\, \fstar)
  \;=\; L\epsilon.
\]
Combining:
\[
  \Pr_{D(\fstar)}(Y \in C_{\hat{f}}(X))
  \;\geq\; 1 - \alpha - L\epsilon.
\]
Equivalently,
\[
  \overline{\mathrm{ExMisc}}(C_{\hat{f}}, D(\fstar)) \leq L\epsilon.
\]
\end{proof}

\subsection{Proof of Parametric Upper Bound (Theorem~\ref{thm:parametric})}

The argument rests on the following predictor-error bound for the
tail-averaged iterate; its proof, a tail-averaging analysis of the
retraining recursion, is given in the supplementary material
(Appendix~\ref{app:polyak}).

\begin{lemma}[Predictor-error bound for P-PCP]\label{lem:polyak-bound}
Under the assumptions of Theorem~\ref{thm:parametric}, the tail-averaged
predictor returned by Algorithm~\ref{alg:pcp} satisfies
\[
  \E\!\left[d_\cF(\bar f_{B,K}, \fstar)^2\right]
  \;\le\;
  \widetilde O\!\left(\frac{1}{(1-\rho)^2\,n}\right).
\]
\end{lemma}

\begin{proof}
The proof has two steps.

\textbf{Step 1: Predictor error after stabilization and averaging.}
By Lemma~\ref{lem:polyak-bound}, Algorithm~\ref{alg:pcp} returns a
predictor $\bar f_{B,K}$ satisfying
\[
  \E\!\left[d_\cF(\bar f_{B,K}, \fstar)^2\right]
  \;\le\;
  \widetilde O\!\left(\frac{1}{(1-\rho)^2\,n}\right).
\]
By Jensen's inequality (since $x \mapsto \sqrt{x}$ is concave),
\[
  \E\!\left[d_\cF(\bar f_{B,K}, \fstar)\right]
  \;\le\;
  \widetilde O\!\left(\frac{1}{(1-\rho)\sqrt n}\right).
\]

\textbf{Step 2: Transfer to coverage.}
Apply Theorem~\ref{thm:meta} with
$\epsilon = d_\cF(\bar f_{B,K}, \fstar)$:
\[
  \overline{\mathrm{ExMisc}}(C_{\bar f_{B,K}}, D(\fstar))
  \;\le\;
  L \cdot \E\!\left[d_\cF(\bar f_{B,K}, \fstar)\right]
  \;\le\;
  \widetilde O\!\left(\frac{L}{(1-\rho)\sqrt n}\right).
\]
Since $L \le 1$ in the regime of interest, this is
\[
  \overline{\mathrm{ExMisc}}(C_{\bar f_{B,K}}, D(\fstar))
  \;\le\;
  \widetilde O\!\left(\frac{1}{(1-\rho)\sqrt n}\right).
\]

When $\rho \asymp L$, for example in the well-conditioned case
$\beta/\gamma \approx 1$, we have
\[
\frac{1}{1-\rho}=\Theta\!\left(\frac{1}{1-L}\right),
\]
so the upper bound matches the $n$- and $L$-dependence of the lower
bound up to logarithmic factors.
\end{proof}

\section{Proof of Model-Free Upper Bound (Theorem~\ref{thm:model-free})}
\label{app:model-free}

We prove the three components: the finite-class contraction argument,
the iterated convergence, and the final coverage guarantee.

\subsection{Finite-Class Contraction Argument}

\begin{proof}
Let $\mathrm{BR}(P) = \arg\min_{f \in \cF} \E_{P}[\ell(f(X),Y)]$
denote the population best response to distribution $P$.  Since
$\fstar$ is performatively stable,
$\fstar = \mathrm{BR}(D(\fstar))$.

\textbf{Step 1: ERM approximation.}
With $m$ i.i.d.\ samples from $D(f_t)$ and $\ell \in [0,1]$,
uniform convergence over $|\cF| = N$ gives:
\[
  \sup_{f \in \cF}
  \bigl|\widehat{R}(f; S_t) - R_{D(f_t)}(f)\bigr|
  \;\leq\; c_0\sqrt{\frac{\log N}{m}}
\]
with probability $\geq 1 - 1/n^2$.  Therefore:
\[
  R_{D(f_t)}(f_{t+1})
  \;\leq\; R_{D(f_t)}(\mathrm{BR}(D(f_t)))
    + 2c_0\sqrt{\frac{\log N}{m}}.
\]

\textbf{Step 2: Best-response stability.}
Define the sub-optimality gap
\[
  \mu \;=\; \min_{f \in \cF,\, f \neq \fstar}
    \bigl[R_{D(\fstar)}(f) - R_{D(\fstar)}(\fstar)\bigr]
  \;>\; 0,
\]
which is strictly positive because $\cF$ is finite and $\fstar$ is
the unique minimizer under $D(\fstar)$.

\emph{Claim:} If $\TV(D(f_t), D(\fstar)) < \mu/2$, then
$\mathrm{BR}(D(f_t)) = \fstar$.

\emph{Proof of claim.}  For any $f \neq \fstar$:
\begin{align*}
  R_{D(f_t)}(f) - R_{D(f_t)}(\fstar)
  &= \bigl[R_{D(\fstar)}(f) - R_{D(\fstar)}(\fstar)\bigr] \\
  &\quad + \bigl[R_{D(f_t)}(f) - R_{D(\fstar)}(f)\bigr]
       - \bigl[R_{D(f_t)}(\fstar) - R_{D(\fstar)}(\fstar)\bigr].
\end{align*}
Since $\ell \in [0,1]$, each difference in the second line is
bounded by $\TV(D(f_t), D(\fstar))$.  Therefore:
\[
  R_{D(f_t)}(f) - R_{D(f_t)}(\fstar)
  \;\geq\; \mu - 2\,\TV(D(f_t), D(\fstar))
  \;>\; 0,
\]
so $\fstar$ is the unique minimizer of $R_{D(f_t)}$ over $\cF$.
\hfill$\square$

By $L$-sensitivity,
$\TV(D(f_t), D(\fstar)) \leq L \cdot d_\cF(f_t, \fstar)$,
so the claim holds whenever
$d_\cF(f_t, \fstar) < \mu/(2L)$.

\textbf{Step 3: Contraction in the near regime.}
When $d_\cF(f_t, \fstar) < \mu/(2L)$, the population best response
is $\fstar$ itself, so
$d_\cF(\mathrm{BR}(D(f_t)), \fstar) = 0$.
Combined with Step~1 and the triangle inequality:
\[
  d_\cF(f_{t+1}, \fstar)
  \;\leq\;
  d_\cF(f_{t+1}, \mathrm{BR}(D(f_t)))
  + d_\cF(\mathrm{BR}(D(f_t)), \fstar)
  \;\leq\;
  O\!\left(\sqrt{\frac{\log N}{m}}\right) + 0.
\]

Since $0 \leq L \cdot d_\cF(f_t, \fstar)$, this implies the stated
recursion:
\[
  d_\cF(f_{t+1}, \fstar)
  \;\leq\; L \cdot d_\cF(f_t, \fstar)
  + O\!\left(\sqrt{\frac{\log N}{m}}\right).
\]

\textbf{Transient.}
The far regime, where $d_\cF(f_t, \fstar) \geq \mu/(2L)$, lasts at most
\[
  T_0
  =
  O\!\left(\frac{\log(\mathrm{diam}(\cF) / \mu)}{1-L}\right)
\]
rounds. This is because the excess risk decreases by at least a
constant per round in that regime.
\end{proof}

\subsection{Iterated Convergence and Coverage}

\begin{proof}[Proof of Theorem~\ref{thm:model-free}]
From the one-step recursion,
\[
d_\cF(f_{t+1},\fstar)
\le
L\,d_\cF(f_t,\fstar)
+
c_1\sqrt{\frac{\log N}{m}},
\]
iterating over $T$ rounds gives
\[
d_\cF(f_T,\fstar)
\le
L^T d_\cF(f_0,\fstar)
+
\sum_{k=0}^{T-1} L^k\,
c_1\sqrt{\frac{\log N}{m}}.
\]
Since $d_\cF(f_0,\fstar)\le \mathrm{diam}(\cF)$ and
$\sum_{k=0}^{T-1} L^k \le 1/(1-L)$, we obtain
\[
d_\cF(f_T,\fstar)
\le
L^T \mathrm{diam}(\cF)
+
\frac{c_1}{1-L}\sqrt{\frac{\log N}{m}}.
\]

\textbf{Choosing $T$.}
Set
\[
T=\left\lceil \frac{\log(n\,\mathrm{diam}(\cF))}{\log(1/L)} \right\rceil
\]
so that
\[
L^T \mathrm{diam}(\cF)\le \frac{1}{n}.
\]
Since $\log(1/L)\ge 1-L$ for $L\in(0,1)$, this choice satisfies
\[
T = O\!\left(\frac{\log n}{1-L}\right).
\]

\textbf{Per-round budget.}
Therefore
\[
m=\left\lfloor \frac{n}{T+1}\right\rfloor
=
\Theta\!\left(\frac{n(1-L)}{\log n}\right).
\]

\textbf{Equilibrium error.}
Substituting this into the previous bound yields
\[
\frac{1}{1-L}\sqrt{\frac{\log N}{m}}
=
\widetilde{O}\!\left(\frac{1}{(1-L)^{3/2}\sqrt n}\right).
\]
Hence
\[
d_\cF(f_T,\fstar)
\le
\frac{1}{n}
+
\widetilde{O}\!\left(\frac{1}{(1-L)^{3/2}\sqrt n}\right)
=
\widetilde{O}\!\left(\frac{1}{(1-L)^{3/2}\sqrt n}\right).
\]

\textbf{Coverage.}
Applying Theorem~\ref{thm:meta} and using the bound above, we conclude
\[
\overline{\mathrm{ExMisc}}(C_{f_T,S},D(\fstar))
\le
\widetilde{O}\!\left(\frac{1}{(1-L)^{3/2}\sqrt n}\right).
\]
\end{proof}

\noindent
The extra $(1-L)^{-1/2}$ factor relative to the minimax lower bound, and
how tail-averaging removes it in the parametric setting, are discussed in
the remark following Theorem~\ref{thm:model-free} in
Section~\ref{sec:model-free}.

\clearpage
\section*{Supplementary Material}
\renewcommand{\theequation}{S\arabic{equation}}
\renewcommand{\thefigure}{S\arabic{figure}}
\renewcommand{\thetable}{S\arabic{table}}
\renewcommand{\thealgorithm}{S\arabic{algorithm}}
\setcounter{equation}{0}
\setcounter{figure}{0}
\setcounter{table}{0}
\setcounter{algorithm}{0}
\renewcommand{\theHequation}{supp.\arabic{equation}}
\renewcommand{\theHfigure}{supp.\arabic{figure}}
\renewcommand{\theHtable}{supp.\arabic{table}}
\renewcommand{\theHalgorithm}{supp.\arabic{algorithm}}
\renewcommand{\theHsection}{supp.\arabic{section}}
\appendix
\renewcommand{\theHsection}{supp.\arabic{section}}
\renewcommand{\thesection}{S\arabic{section}}

\section{Tail-Averaging Analysis for the Parametric Algorithm}
\label{app:polyak}

This appendix proves Lemma~\ref{lem:polyak-bound}, the predictor-error
bound used in the proof of Theorem~\ref{thm:parametric}.

The main idea is the following.
Repeated retraining moves the predictor toward equilibrium because the
population retraining map is contractive.
However, each empirical retraining step is noisy because it uses only
finitely many samples.
To combine these two facts, we split the procedure into two phases:

\begin{enumerate}[leftmargin=*,itemsep=2pt]
\item a short \emph{stabilization period}, whose only role is to reduce
the effect of the initial predictor, and
\item an \emph{averaging period}, where we average later iterates that
are already close to equilibrium in order to reduce sampling noise.
\end{enumerate}

The appendix below makes this idea precise.

\subsection{Parameter notation and one-step decomposition}
\label{app:polyak-linearization}

For the parametric analysis, write the predictor as $f_\theta$ with
parameter $\theta \in \Theta \subseteq \R^d$, and take
\[
d_\cF(f_\theta, f_{\theta'}) = \|\theta-\theta'\|_2.
\]
Write $\theta_t$ for the parameter of the iterate $f_t$, and let
$\theta^\star$ denote the equilibrium parameter corresponding to
$\fstar$.
Define the error
\[
e_t \coloneqq \theta_t - \theta^\star.
\]

Let
\[
\mathrm{BR}(D(\theta))
\;\coloneqq\;
\arg\min_{\vartheta \in \Theta} R_{D(\theta)}(\vartheta)
\]
be the population best response to deployment $\theta$.
Under $\beta$-smoothness, $\gamma$-strong convexity, and
$L$-sensitivity, the best-response map is contractive:
\begin{equation}
\label{eq:app-br-contraction}
\|\mathrm{BR}(D(\theta))-\theta^\star\|
\;\le\;
\rho \|\theta-\theta^\star\|,
\qquad
\rho = \frac{L\beta}{\gamma} < 1,
\end{equation}
by \citet[Theorem~3.5]{perdomo2020performative}.

At round $t$, Algorithm~\ref{alg:pcp} computes an ERM from $m$ samples
drawn under $D(\theta_t)$.
Under standard finite-sample regularity for ERM, we may write
\begin{equation}
\label{eq:app-erm-decomp}
\theta_{t+1}
=
\mathrm{BR}(D(\theta_t)) + \xi_t,
\end{equation} 
where the estimation error $\xi_t$ satisfies
\begin{equation}
\label{eq:app-erm-noise}
\E[\xi_t \mid \theta_t] = 0,
\qquad
\E[\|\xi_t\|^2 \mid \theta_t] \le \frac{\sigma^2}{m}
\end{equation}
for a constant $\sigma^2$ that depends on the parametric problem but
not on $t$, $m$, or $n$.

Equation~\eqref{eq:app-erm-decomp} says that each empirical retraining
step equals the population retraining step plus a finite-sample error.

\begin{lemma}[Contractive Error Recursion]
\label{lem:linearization}
Under \eqref{eq:app-br-contraction}--\eqref{eq:app-erm-noise}, the
errors satisfy
\[
e_{t+1} = A_t e_t + \xi_t,
\]
where $\|A_t\| \le \rho$ for all $t$.
Consequently, for every $t > s$ there are matrices
$\Phi_{t,s}$ with
\[
\|\Phi_{t,s}\| \le \rho^{\,t-s}
\]
such that
\[
e_t
=
\Phi_{t,0} e_0
+
\sum_{s=0}^{t-1} \Phi_{t,s+1}\,\xi_s.
\]
\end{lemma}

\begin{proof}
From \eqref{eq:app-erm-decomp},
\[
e_{t+1}
=
\theta_{t+1}-\theta^\star
=
\bigl(\mathrm{BR}(D(\theta_t))-\theta^\star\bigr) + \xi_t.
\]
Let
\[
u_t \coloneqq \mathrm{BR}(D(\theta_t))-\theta^\star.
\]
By \eqref{eq:app-br-contraction},
\[
\|u_t\| \le \rho \|e_t\|.
\]

If $e_t\neq 0$, define
\[
A_t \coloneqq \frac{u_t e_t^\top}{\|e_t\|^2}.
\]
Then
\[
A_t e_t = u_t,
\qquad
\|A_t\| = \frac{\|u_t\|}{\|e_t\|} \le \rho.
\]
If $e_t=0$, then $u_t=0$ by the contraction bound, so we may set
$A_t=0$.
Thus in all cases
\[
u_t = A_t e_t
\qquad\text{with}\qquad
\|A_t\|\le \rho.
\]
Therefore
\[
e_{t+1} = A_t e_t + \xi_t,
\]
which proves the first claim.

Iterating the recursion yields
\[
e_t
=
\Phi_{t,0} e_0
+
\sum_{s=0}^{t-1} \Phi_{t,s+1}\,\xi_s,
\]
where
\[
\Phi_{t,s}
=
A_{t-1}A_{t-2}\cdots A_s
\qquad (t>s),
\]
and $\Phi_{t,t}=I$.
Since each $\|A_j\|\le \rho$,
\[
\|\Phi_{t,s}\|
\le
\prod_{j=s}^{t-1}\|A_j\|
\le
\rho^{\,t-s}.
\]
\end{proof}

\subsection{Error of the averaged post-stabilization iterates}
\label{app:polyak-variance}

Algorithm~\ref{alg:pcp} uses two phases:

\begin{itemize}[itemsep=2pt]
\item the first $B$ rounds reduce the effect of the initial predictor;
\item the next $K$ rounds are averaged.
\end{itemize}

Define the averaged parameter error over the second phase by
\[
\bar e_{B,K}
\;\coloneqq\;
\frac{1}{K}\sum_{t=B+1}^{B+K} e_t.
\]

The next lemma separates the two sources of error:

\begin{enumerate}[leftmargin=*,itemsep=2pt]
\item an \emph{initialization term}, which becomes small after the
stabilization period, and
\item a \emph{sampling-noise term}, which becomes small after averaging.
\end{enumerate}

\begin{lemma}[Bias-Variance Bound for the Post-Stabilization Average]
\label{lem:averaged-variance}
Under the conditions of Lemma~\ref{lem:linearization},
\[
\E\!\left[\|\bar e_{B,K}\|^2\right]
\;\le\;
2\rho^{2B}\|e_0\|^2
+
\frac{2\sigma^2(B+K)}{mK^2(1-\rho)^2}.
\]
\end{lemma}

\begin{proof}
Using Lemma~\ref{lem:linearization},
\[
\bar e_{B,K}
=
\frac{1}{K}\sum_{t=B+1}^{B+K}\Phi_{t,0}e_0
+
\frac{1}{K}\sum_{t=B+1}^{B+K}\sum_{s=0}^{t-1}\Phi_{t,s+1}\xi_s.
\]
Call these two terms $\bar b_{B,K}$ and $\bar v_{B,K}$:
\[
\bar e_{B,K} = \bar b_{B,K} + \bar v_{B,K}.
\]

\paragraph{Initialization term.}
For every $t \ge B+1$,
\[
\|\Phi_{t,0}e_0\|
\le
\rho^t \|e_0\|
\le
\rho^B \|e_0\|.
\]
Therefore
\[
\|\bar b_{B,K}\|
=
\left\|
\frac{1}{K}\sum_{t=B+1}^{B+K}\Phi_{t,0}e_0
\right\|
\le
\frac{1}{K}\sum_{t=B+1}^{B+K}\rho^B\|e_0\|
=
\rho^B\|e_0\|.
\]
Hence
\[
\|\bar b_{B,K}\|^2 \le \rho^{2B}\|e_0\|^2.
\]

\paragraph{Sampling-noise term.}
Rearrange the second term by grouping together the contributions of
each noise variable $\xi_s$:
\[
\bar v_{B,K}
=
\sum_{s=0}^{B+K-1} W_s \xi_s,
\]
where
\[
W_s
\;\coloneqq\;
\frac{1}{K}
\sum_{t=\max(B+1,s+1)}^{B+K}\Phi_{t,s+1}.
\]
Using $\|\Phi_{t,s+1}\|\le \rho^{\,t-s-1}$,
\[
\|W_s\|
\le
\frac{1}{K}
\sum_{j=0}^{\infty}\rho^j
=
\frac{1}{K(1-\rho)}.
\]
Condition on the full trajectory $\mathcal{A} = \{A_0, \ldots, A_{B+K-1}\}$.
Given $\mathcal{A}$, each $W_s$ is deterministic.
Since the samples at each round are fresh and independent,
$\E[\xi_s \mid \mathcal{A}] = 0$ and the $\xi_s$ are
conditionally independent given $\mathcal{A}$, so the cross-terms vanish:
\[
\E\!\left[\|\bar v_{B,K}\|^2 \mid \mathcal{A}\right]
= \sum_{s=0}^{B+K-1} \|W_s\|^2\,\E\!\left[\|\xi_s\|^2 \mid \mathcal{A}\right]
\le \frac{\sigma^2}{m}\sum_{s=0}^{B+K-1} \|W_s\|^2.
\]
Taking expectation over $\mathcal{A}$ and applying
$\|W_s\| \le \frac{1}{K(1-\rho)}$ gives
\[
\E\!\left[\|\bar v_{B,K}\|^2\right]
\le
(B+K)\cdot \frac{1}{K^2(1-\rho)^2}\cdot \frac{\sigma^2}{m}.
\]

\paragraph{Combine the two terms.}
Using $\|a+b\|^2 \le 2\|a\|^2 + 2\|b\|^2$,
\[
\E\!\left[\|\bar e_{B,K}\|^2\right]
\le
2\|\bar b_{B,K}\|^2 + 2\E\!\left[\|\bar v_{B,K}\|^2\right],
\]
which gives
\[
\E\!\left[\|\bar e_{B,K}\|^2\right]
\le
2\rho^{2B}\|e_0\|^2
+
\frac{2\sigma^2(B+K)}{mK^2(1-\rho)^2}.
\]
\end{proof}

\begin{remark}[Why the Stabilization Period Is Needed]
If we averaged \emph{all} iterates, then the average would include the
early rounds that are still strongly affected by the initial predictor.
The stabilization period removes these early terms.
After that, averaging helps because it reduces the noise
of later iterates that are already close to equilibrium.
\end{remark}

\subsection{Rate for Algorithm~\ref{alg:pcp}}
\label{app:polyak-combining}

The last step is to choose $B$ and $K$ so that:

\begin{itemize}[itemsep=2pt]
\item the initialization term is negligible, and
\item the averaging term uses essentially the full sample budget.
\end{itemize}

\begin{lemma}[Tail-Average Rate for Performative Retraining]
\label{lem:polyak-ruppert}
Under the conditions of Theorem~\ref{thm:parametric}, suppose the ERM
errors satisfy \eqref{eq:app-erm-noise}.  Set
\[
B = K = \left\lceil \frac{\log n}{\log(1/\rho)} \right\rceil,
\qquad
m = \left\lfloor \frac{n}{B+K+1} \right\rfloor.
\]
Then the predictor returned by Algorithm~\ref{alg:pcp} satisfies
\[
\E\!\left[d_\cF(\bar f_{B,K}, \fstar)^2\right]
\;\le\;
\widetilde O\!\left(\frac{1}{(1-\rho)^2\,n}\right).
\]
Consequently,
\[
\E\!\left[d_\cF(\bar f_{B,K}, \fstar)\right]
\;\le\;
\widetilde O\!\left(\frac{1}{(1-\rho)\sqrt n}\right).
\]
\end{lemma}

\begin{proof}
Apply Lemma~\ref{lem:averaged-variance}.

\paragraph{Initialization term.}
With
\[
B = \left\lceil \frac{\log n}{\log(1/\rho)} \right\rceil,
\]
we have
\[
\rho^B \le \frac{1}{n},
\qquad\text{hence}\qquad
\rho^{2B}\|e_0\|^2 \le \frac{\|e_0\|^2}{n^2}.
\]
So the initialization term is of order $1/n^2$, which is negligible
compared to the sampling-noise term of order $1/n$ that dominates the
final bound.

\paragraph{Averaging term.}
Since $B=K$, the factor $(B+K)/K^2$ is of order $1/K$.
Also, because
\[
m = \left\lfloor \frac{n}{B+K+1} \right\rfloor,
\]
we have $m \asymp n/(B+K)$, so
\[
\frac{B+K}{mK^2}
=
O\!\left(\frac{1}{n}\right).
\]
Therefore
\[
\frac{2\sigma^2(B+K)}{mK^2(1-\rho)^2}
=
O\!\left(\frac{1}{n(1-\rho)^2}\right).
\]

\paragraph{Combine.}
Lemma~\ref{lem:averaged-variance} now gives
\[
\E\!\left[\|\bar e_{B,K}\|^2\right]
\le
O\!\left(\frac{1}{n^2}\right)
+
O\!\left(\frac{1}{n(1-\rho)^2}\right)
=
\widetilde O\!\left(\frac{1}{n(1-\rho)^2}\right).
\]
Since $d_\cF$ is the Euclidean distance in parameter space for the
parametric analysis, this is 
\[
\E\!\left[d_\cF(\bar f_{B,K}, \fstar)^2\right]
\le
\widetilde O\!\left(\frac{1}{n(1-\rho)^2}\right).
\]
The bound on the first moment follows from Jensen's inequality:
\[
\E\!\left[d_\cF(\bar f_{B,K}, \fstar)\right]
\le
\sqrt{
\E\!\left[d_\cF(\bar f_{B,K}, \fstar)^2\right]
}
\le
\widetilde O\!\left(\frac{1}{(1-\rho)\sqrt n}\right).
\]
\end{proof}

\begin{remark}[Interpretation of the Choice of $B$ and $K$]
The stabilization length $B$ is chosen so that the residual effect of
the initial predictor, which decays geometrically like $\rho^B$, is
smaller than the statistical error level of the final estimator.

The averaging window $K$ is chosen to be of the same logarithmic order.
This gives enough post-stabilization iterates to reduce sampling noise
through averaging, while still keeping the per-round sample size
\[
m=\lfloor n/(B+K+1)\rfloor
\]
large enough for accurate ERM at each round.

Thus both $B$ and $K$ need only be of order $\log n$: geometric
contraction makes the initialization bias negligible after
$O(\log n)$ rounds, and a logarithmic averaging window is sufficient
to control the remaining noise.
\end{remark}

\section{Experimental Details}
\label{app:experiments}

This appendix collects the full settings for the simulations reported in
Section~\ref{sec:experiments}.  Unless stated otherwise, all
simulations use $\alpha=0.1$ (target coverage $90\%$) and
50{,}000 independent test samples drawn from the equilibrium
distribution $D(\fstar)$.

For the semi-synthetic strategic-regression setup, we also align with the broader strategic-classification and recourse literature on behavioral response and practical deployment constraints~\citep{hu2019disparate,lechner2023strategic,ustun2019actionable}.

\subsection{Parameter Specifications}
\label{app:params}

Table~\ref{tab:params} summarizes the parameter grids for the six
simulation setups reported in the main text and appendix.  All use a constant
covariate, $\cX=\{x_0\}$, so that the simulations isolate the
performative effect without additional variation from features.

\begin{table}[ht]
\caption{Synthetic experiments reported in the paper and appendix.}
\label{tab:params}
\centering
\small
\begin{tabular}{@{}llp{5.4cm}c@{}}
\toprule
\textbf{Simulation} & \textbf{Figures} & \textbf{Purpose / key parameters} & \textbf{Reps} \\
\midrule
1 & Fig.~\ref{fig:synthetic_all}(a,b) &
$L \in \{0,0.1,\ldots,0.7\}$,
$n \in \{50,100,200,500,1000,3000,10^4,3\times 10^4,10^5\}$,
Gaussian location model, $a=0.3$;
tests rate and price of performativity
& 2000 \\[3pt]

2 & Fig.~\ref{fig:synthetic_all}(c) &
$L \in \{0.1,0.3,0.5,0.7,0.9\}$,
$\varepsilon \in \{0,0.05,0.1,0.2,0.5,1,2\}$,
$m=2000$, Gaussian location model, $a=0.3$;
tests $L\varepsilon$ collapse
& 2000 \\[3pt]

3 & Fig.~\ref{fig:phase_transition} &
\parbox[t]{5.4cm}{\raggedright
$L \in \{0.7,0.8,0.9,0.95,0.98,0.99,$\\
$1.0,1.01,1.02,1.05,1.1,1.2\}$,\\
bounded Bernoulli sweep at $n=300{,}000$ with\\
$c=0.05$ and $\varepsilon=0.01$;\\
tests the phase boundary with explicit equilibria for all $L$}
& 2000 \\[3pt]

4 & Fig.~\ref{fig:sim1} &
$L \in \{0,0.3,0.5,0.7,0.9\}$,
$n \in \{100,500,1000,5000,10^4,5\times 10^4\}$,
$\delta=0.5/\sqrt{n}$; Le~Cam two-hypothesis construction;
worst-case excess miscoverage and scaled rate check
& 2000 \\[3pt]

5 & Fig.~\ref{fig:sim3} &
Finite-class Gaussian location model for MF-PCP;
$L \in \{0,0.3,0.5,0.7\}$,
$N \in \{10,50,200\}$,
$n \in \{200,500,1000,5000,10^4,5\times 10^4\}$
& 2000 \\[3pt]

6 & Fig.~\ref{fig:sim6} &
Impossibility construction at $L=1$ using $D_+$ and $D_-$,
$n \in \{100,500,1000,5000,10^4,5\times 10^4,10^5\}$,
$\varepsilon=0.5/n$
& 2000 \\
\bottomrule
\end{tabular}
\end{table}

\paragraph{Distribution map.}
Simulations~1, 2, 4, and~5 use the Gaussian location model
\[
Y \sim \cN(a + L\theta, 1),
\]
with equilibrium
\[
\theta^\star = \frac{a}{1-L},
\]
and score
\[
s(y;\theta)=|y-\theta|.
\]

Simulation~1 varies $L$ and $n$ to study the rate
$\widetilde{\Theta}(1/((1-L)\sqrt{n}))$ and the induced sample
complexity scaling $\widetilde{\Theta}(1/(1-L)^2)$ for a fixed target
accuracy.

Simulation~2 fixes the calibration size $m=2000$ and varies the
deployed predictor offset
\[
\varepsilon = |\theta_{\mathrm{deployed}}-\theta^\star|
\]
across several values of $L$ in order to test the transfer-theorem prediction that small coverage losses are mainly
determined by the product $L\varepsilon$.

\section{Auxiliary Experiments}\label{app:additional}

\subsection{Simulation~4: Le Cam Lower-Bound}

This simulation validates the minimax lower bound of
Theorem~\ref{thm:minimax} by instantiating the Le~Cam
two-hypothesis construction.  For each $(L, n)$, we run
Algorithm~\ref{alg:pcp} separately on maps $D_0$ and $D_1$
(differing by $\delta = 0.5/\!\sqrt{n}$) and report the
worst-case (maximum) excess miscoverage.

\begin{figure}[tbp]
    \centering

    \begin{subfigure}[t]{0.49\linewidth}
        \centering
        \includegraphics[width=\linewidth]{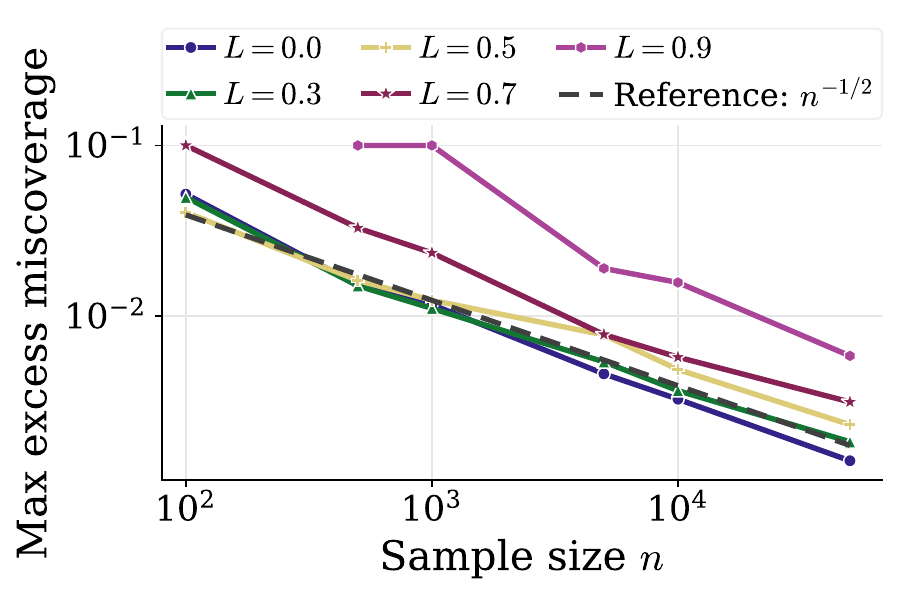}
        \vspace{-1.5em}
        \caption{Worst-case excess miscoverage versus $n$}
        \label{fig:sim1_worstcase}
    \end{subfigure}\hfill
    \begin{subfigure}[t]{0.49\linewidth}
        \centering
        \includegraphics[width=\linewidth]{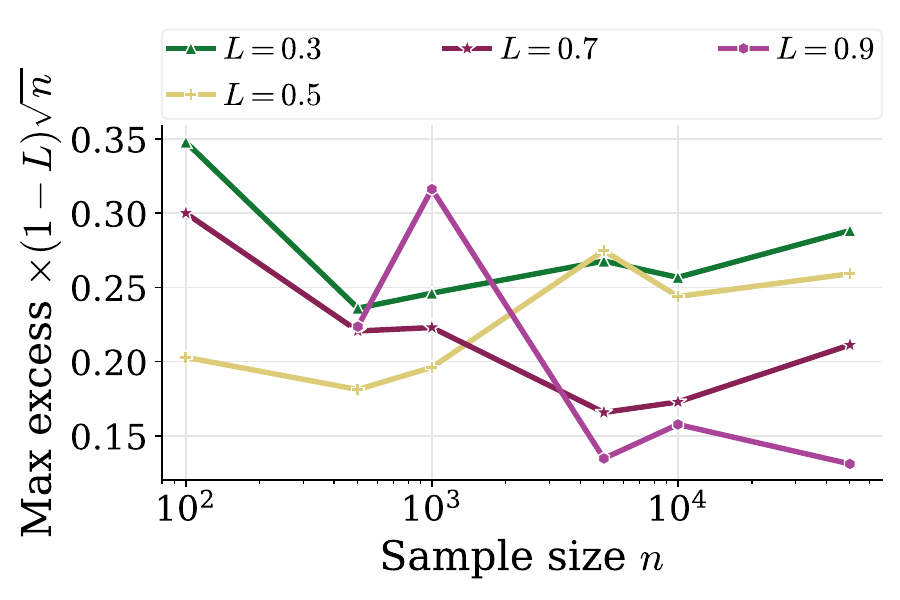}
        \vspace{-1.5em}
        \caption{Scaled excess miscoverage}
        \label{fig:sim1_scaled}
    \end{subfigure}

    \caption{Simulation~4: Le~Cam lower bound validation.
    (a) Worst-case excess miscoverage (maximum over $D_0, D_1$)
    vs.\ $n$ on a log-log scale.  The $L = 0.9$ curve saturates at
    the trivial level for $n \leq 1000$ (algorithm has not converged).
    (b) Scaled quantity
    $\mathrm{excess} \times (1{-}L) \times \sqrt{n}$; approximate
    flatness in $[0.13, 0.35]$ confirms the
    $\Omega(1/((1{-}L)\sqrt{n}))$ rate.}
    \label{fig:sim1}
\end{figure}

Figure~\ref{fig:sim1}\subref{fig:sim1_worstcase} shows the worst-case excess miscoverage
on a log-log scale.  All curves decrease with $n$, with higher $L$
producing larger miscoverage at each sample size.  The $L = 0.9$
curve starts at the trivial level (excess $= \alpha$) for
$n \leq 1000$ because the algorithm has not converged, but
recovers the correct rate for $n \geq 5000$.

Figure~\ref{fig:sim1}\subref{fig:sim1_scaled} plots the scaled quantity
$\text{excess} \times (1-L) \times \sqrt{n}$, which should be
approximately constant if the $\Omega(1/((1{-}L)\sqrt{n}))$ rate
is tight.  The scaled values remain in $[0.13, 0.35]$ across
three orders of magnitude in $n$ for all $L > 0$, confirming the
rate.  The scaled constant decreases mildly with $L$
(from ${\approx}\,0.27$ at $L = 0.3$ to ${\approx}\,0.14$ at
$L = 0.9$), which reflects the fact that the Gaussian/Bernoulli
model is not the adversarial worst case (the proof's lower bound
construction uses a different mechanism), but the \emph{rate} in
both $n$ and $(1-L)$ is confirmed.

\subsection{Simulation~5: Model-Free Algorithm}

This simulation validates Algorithm~\ref{alg:mfpcp} (MF-PCP) and
Theorem~\ref{thm:model-free} using the Gaussian location model with a finite
predictor class of size $N$.

\begin{figure}[tbp]
    \centering

    \begin{subfigure}[t]{0.49\linewidth}
        \centering
        \includegraphics[width=\linewidth]{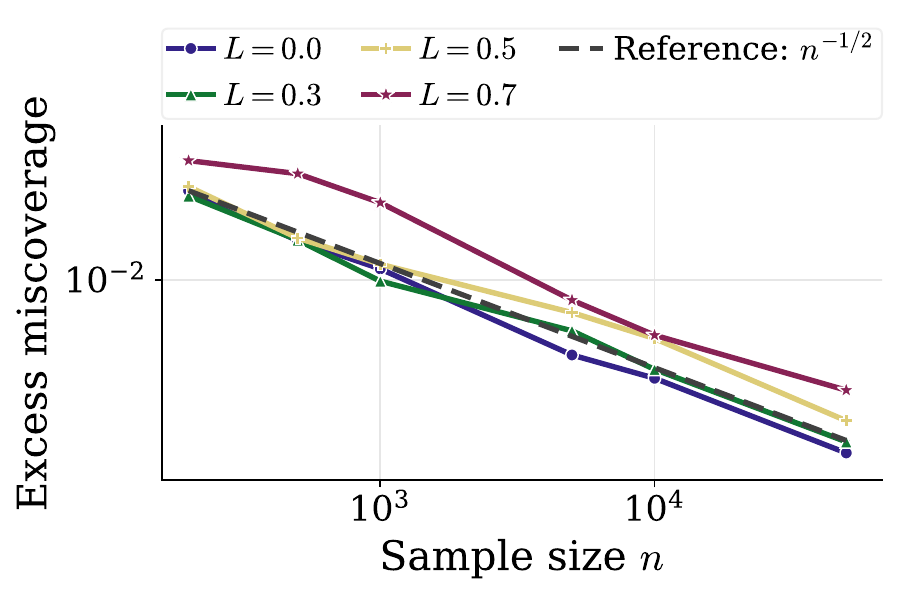}
        \vspace{-1.5em}
        \caption{Excess miscoverage versus $n$ for $N = 50$}
        \label{fig:sim3_convergence}
    \end{subfigure}\hfill
    \begin{subfigure}[t]{0.49\linewidth}
        \centering
        \includegraphics[width=\linewidth]{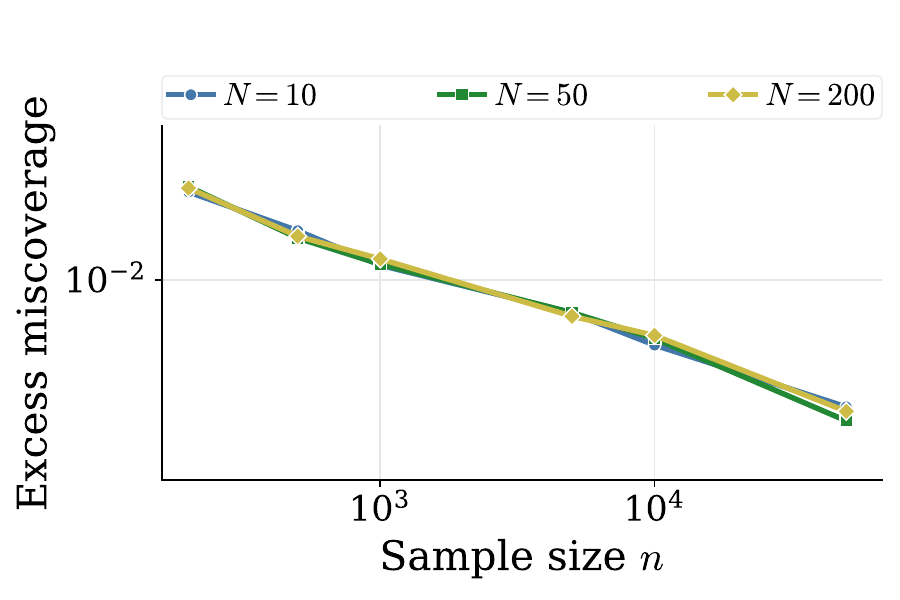}
        \vspace{-1.5em}
        \caption{Effect of class size $N$ at $L = 0.5$}
        \label{fig:sim3_classsize}
    \end{subfigure}

    \caption{Simulation~5: Model-free algorithm validation.
    (a) Excess miscoverage vs.\ $n$ for $N = 50$, confirming
    $n^{-1/2}$ convergence across $L$ values.
    (b) Effect of class size $N$ at $L = 0.5$; curves nearly
    overlap, indicating the $\log N$ factor is small.}
    \label{fig:sim3}
\end{figure}

Figure~\ref{fig:sim3}\subref{fig:sim3_convergence} shows excess miscoverage vs.\ $n$ for
$N = 50$ across four values of $L$.  All curves decrease at
approximately the $n^{-1/2}$ rate (log-log slopes in
$[-0.54, -0.44]$), consistent with the theorem.  Higher $L$ shifts
curves upward, as expected.

Figure~\ref{fig:sim3}\subref{fig:sim3_classsize} shows the effect of the class size $N$
at fixed $L = 0.5$.  The three curves ($N = 10, 50, 200$) are
nearly indistinguishable, confirming that the $\log N$ factor in
the bound is small relative to the dominant $n^{-1/2}$ term for
the class sizes considered.

\paragraph{Convergence to equilibrium.}
The convergence check (fraction of replications where
$\theta_{\text{final}}$ equals the equilibrium grid point
$\theta^\star_{\text{grid}}$) is low across all configurations
(maximum 51\% at $L = 0$, $N = 10$, $n = 50{,}000$; as low as
2.5\% at $L = 0.5$, $N = 200$, $n = 200$).  This does \emph{not}
indicate a problem: with a discrete grid of $N$ points and
stochastic ERM, the algorithm frequently selects a neighbor of
$\theta^\star_{\text{grid}}$ rather than the exact grid point.
Crucially, coverage validity does not require exact convergence; it
only requires that the predictor is close enough that the conformal
calibration set is representative, which is confirmed by the
decreasing excess miscoverage.

\begin{figure}[tbp]
    \centering
    \includegraphics[width=0.50\columnwidth]{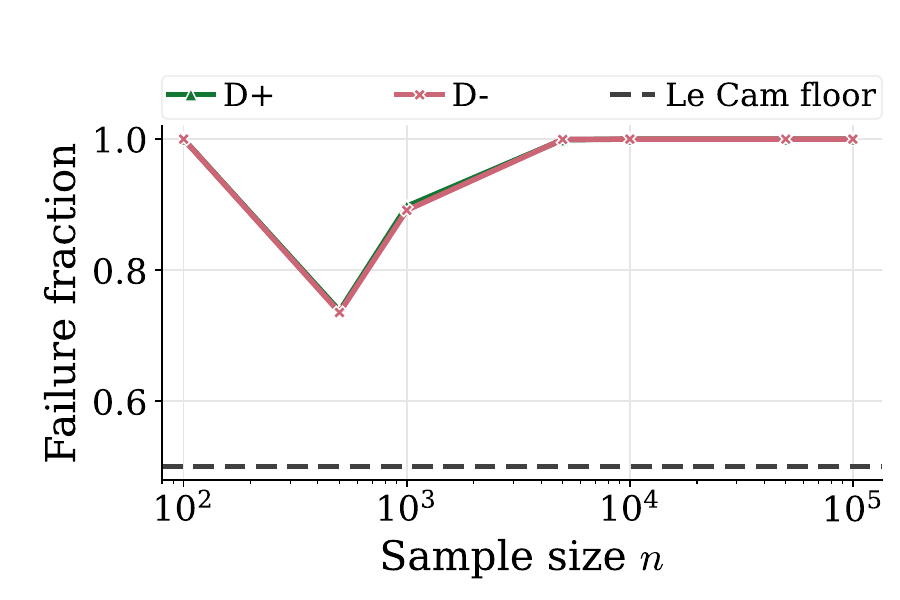}
    \caption{Simulation~6: Impossibility at $L = 1$.
    Failure rate (trivial set or coverage failure) vs.\ $n$ for
    maps $D_+$ and $D_-$. The rate stays bounded away from zero
    at all sample sizes, consistent with Theorem~\ref{thm:impossibility}.
    Dashed line: Le Cam floor at $0.5$.}
    \label{fig:sim6}
\end{figure}

\subsection{Simulation~6: Impossibility at \texorpdfstring{$L=1$}{L=1}}

This simulation validates Theorem~\ref{thm:impossibility} (Case~2)
using the construction from the proof.  Two 1-sensitive maps
$D_+$ and $D_-$ have equilibria at $\theta^\star_+ = 1 - \varepsilon$
and $\theta^\star_- = \varepsilon$, respectively, where
$\varepsilon = 0.5/n$.  As $n$ grows, the two maps become
increasingly similar (since $\varepsilon \to 0$), yet their
equilibria remain far apart ($\theta^\star_+ \to 1$,
$\theta^\star_- \to 0$).  The impossibility theorem predicts that
no algorithm can simultaneously achieve non-trivial prediction sets
and valid coverage under both maps.

Figure~\ref{fig:sim6} plots the failure rate (fraction of
replications producing either a trivial set $\{0,1\}$ or
coverage below $1 - \alpha$) against $n$.  The failure rate
is 100\% at $n = 100$ and $n \geq 5{,}000$ for both maps.
At intermediate sample sizes ($n = 500, 1000$), the rate dips
to ${\approx}\,74\%$ and ${\approx}\,89\%$ respectively, because
with larger $\varepsilon$ the algorithm occasionally identifies
the correct map and produces a valid non-trivial set.  As $n$
increases further, $\varepsilon \to 0$ makes the maps
indistinguishable and the failure rate returns to 100\%. 

This non-monotone pattern is predicted by the construction:
at small $n$, $\varepsilon = 0.5/n$ is large enough that
$\TV(D_+(\theta), D_-(\theta)) \leq 2\varepsilon$ remains
detectable, so the algorithm occasionally identifies the
correct map; as $n$ grows, $\varepsilon \to 0$ drives the
total KL below the discrimination threshold, restoring
impossibility.

The dominant failure mode is trivial prediction sets: at
$n \geq 5{,}000$, the fraction producing $C = \{0,1\}$ exceeds
99.7\%.  This confirms the theorem's prediction that the
algorithm hedges by outputting the full label space when it
cannot distinguish the maps.

\end{document}